\documentclass[11pt,a4paper]{article}
\usepackage{amsfonts}
\usepackage{amsmath,amssymb}
\usepackage{mathrsfs}
\usepackage{hyperref}
\usepackage{graphicx}
\usepackage{float}
\usepackage{indentfirst}
\newtheorem{thm}{Theorem}[section]

\newtheorem{lem}[thm]{Lemma}
\newtheorem{prop}[thm]{Proposition}
\newtheorem{defn}[thm]{Definition}

\newtheorem{rem}[thm]{Remark}

\numberwithin{equation}{section}\allowdisplaybreaks

\def\leq{\leqslant}
\def\ge{\geqslant}

\def\leq{\leqslant}
\def\geq{\geqslant}

\def\Z{{\mathbb{Z}}}

\begin{document}

\title{\large\bf  Uniform local  well-posedness and inviscid limit for the Benjamin-Ono-Burgers equation}

\author{\normalsize \bf Mingjuan Chen$^a$,\quad Boling Guo$^a$, \quad Lijia Han$^{b,}$\footnote{Corresponding author.} \\
\footnotesize
\it a. Institute of Applied Physics and Computational Mathematics, Beijing 100088, PR China; \\
\footnotesize
\it b. Department of Mathematics and Physics, North China Electric Power University, Beijing 102206,  PR China.\\
\footnotesize
\it Emails: mjchenhappy@pku.edu.cn; gbl@iapcm.ac.cn; hljmath@ncepu.edu.cn \ \  \\
}
\date{} \maketitle
\begin{minipage}{13.5cm}
\footnotesize \bf Abstract. \rm  In this paper, we study the Cauchy
problem for the Benjamin-Ono-Burgers equation $\partial_t u-\epsilon \partial_x^2 u+\mathcal{H}\partial_x^2u+u u_x=0$, where $\mathcal{H}$ denotes the Hilbert
transform. We obtain that
it is uniformly locally well-posed for small data in the refined Sobolev space $\widetilde{H}^\sigma(\mathbb{R})$($\sigma\geq 0$), whose low-frequency part is scaling critical and high-frequency part is equal to Sobolev space $H^\sigma$($\sigma\geq 0$). Furthermore, we also obtain its inviscid limit behavior in $\widetilde{H}^\sigma(\mathbb{R})$($\sigma\geq 0$).
\\

\\

\end{minipage}

\section{Introduction}\label{section1}

In this paper, we study the Cauchy problem for the Benjamin-Ono-Burgers (BOB) equation on the real line
\begin{equation}\label{BOB}
\begin{cases}
\partial_t u-\epsilon \partial_x^2 u+\mathcal{H}\partial_x^2u+u u_x=0,\ \ (x,t)\in\mathbb{R}\times\mathbb{R}^+,\\
u(x,0)=\phi(x),
\end{cases}
\end{equation}
where $0<\epsilon\leq 1$, $u$ is a real-valued function of
$(x, t) \in \mathbb{R}\times \mathbb{R}^+$, $\mathcal{H}$ is the Hilbert transform operator defined as follows
\begin{equation*}
\mathcal{H}(f)(x)=p.v.\frac{1}{\pi}\int_\mathbb{R}\frac{f(y)}{x-y}dy.
\end{equation*}

When $\epsilon= 0$, the equation \eqref{BOB} reduces to the classical Benjamin-Ono(BO) equation
\begin{align}\label{BO}
\partial _{t}u+\mathcal{H}\partial_x^2 u+u u_x=0, \quad
u(x,0)=  \phi(x),
\end{align}
which was originally derived as a model in
the study of one-dimensional long internal gravity waves in deep stratified fluids with great depth $\cite{BTB,
Ono}$. The BOB model \eqref{BOB} was obtained by Ewdin and Roberts \cite{ER} in the study of intense magnetic flux tubes of the solar atmosphere. The dissipative effects $-\epsilon \partial_x^2 u$ in that literature are due to weak thermal conduction, where $\epsilon$ is a measure of
the importance of thermal conduction and is assumed small.

Recently, there are many authors who devoted themselves to studying the well-posedness theory and limit behavior of BO and BOB equations. The best result so far  for global well-posedness  of BO equation was proved by Ionescu and Kenig \cite{Io-Ke}  in Sobolev space $H^\sigma$, $\sigma\geq 0$.
For BOB equation, thanks to the dissipative effects, there are many results about its wellposedness. Otani \cite{MOT} derived the global
well-posedness in $H^\sigma$ for $\sigma>-1/2$ by using the Picard methods.  Vento \cite{Vento} proved this result is critical in the sense that the mapping data-solution fails to be $C^3$ continuous if $\sigma<-1/2$.
For more results, we refer to \cite{ABFS,Iorio,KK,KPV,Koch0,Koch,MRB,MRB0,MRB1,molinet3,Ponce,Tao} and the references therein. However, if we consider the uniform well-posedness and inviscid limit for the solutions
of BOB equation, the dissipative effects which related to $\epsilon$ could not be used.

In \cite{Tao}, Tao conjectured  it is feasible to prove that the solutions of BOB
equation converge to those of BO equation when $\epsilon \rightarrow 0$ .  Motivated by \cite{Io-Ke} and \cite{Tao}, Guo and his co-authors \cite{GLB} obtained that BOB equations were uniformly globally well-posed in $H^\sigma$ for $\sigma\geq 1$ and the solutions of BOB converged
to those of BO in $C([0, T]: H^\sigma) (\sigma\geq 1$) for any $T >0$. This result was improved to the energy space $H^{1/2}$ by Molinet \cite{molinet2}. In the light of \cite{Io-Ke}, it seems natural to obtain the limit behavior of the real-valued solutions to BOB equation in $H^\sigma$, $\sigma\geq 0$. To the best of our knowledge, the limit behavior of BOB equation in $H^\sigma$($0 \leq \sigma < \frac{1}{2}$) is still open. Our main goal in this paper is to fill the gap between $L^2$ and $H^{1/2}$.

We obtain that BOB equation is uniformly locally well-posed for small data in the refined Sobolev space $\widetilde{H}^\sigma$($\sigma\geq 0$),  whose low-frequency part is scaling critical and high-frequency part is equal to Sobolev space $H^\sigma$($\sigma\geq 0$). In fact, the high-frequency part has already reduced to $L^2$, while the low-frequency part has some special structure. For BO equation, the special structure can be eliminated by performing a gauge transformation in \cite{Io-Ke}. However, this gauge transformation is not available for BOB equation, due to the dissipative structure. We notice that both \cite{GLB} and \cite{molinet2} did not apply gauge transformation.

The basic ideas for the inviscid limit are to get the uniform well-posedness and difference estimates. We first use similar spaces as that in \cite{Io-Ke} which considered BO equation to obtain the bilinear estimates. In order to weaken the interaction between very low and very high frequencies, which is out of control by standard Bourgain method, we assume that low-frequency functions have some additional structure(see the definitions of $X_0$,$Y_0$ and $B_0$). To avoid the logarithmic divergences we work with high-frequency functions that have two components: a weighted $X^{\sigma,b}$-type component(see $X_k$) and a normalized $L^1_xL^2_t$ component(see $Y_k$) which related to smoothing effect. This type of spaces have been used in \cite{Io-Ke,Io-Ke2} and the references about wave maps therein.

Different from \cite{Io-Ke}, we have to construct the uniform homogeneous and inhomogeneous linear estimates for BOB equation. The dissipative structure destroys some symmetries and brings some logarithmic divergences, which will bring several technical difficulties to obtain the uniform estimates. In order to avoid the logarithmic divergence, the homogeneous dyadic decomposition is performed to construct the low-frequency space $Y_0$. Specifically, we need to conquer the singularity which occurs in low-frequency low-modulation cases, when treating  $1/(\tau-\omega(\xi)-i\epsilon\xi^2)$. We lead the readers to Lemma \ref{lemma3} and the proof of Lemma \ref{Lemma2}.  We believe that these techniques can be used in some other problems.

 Let $\mathcal{F}$($\mathcal{F}^{-1}$) denote the (inverse) Fourier transform operators on $\mathcal{S}'(\mathbb{R}\times\mathbb{R})$. Let $\mathcal{F}_x$($\mathcal{F}_\xi^{-1}$) and $\mathcal{F}_t$($\mathcal{F}_\tau^{-1}$) denote the (inverse) Fourier transform operators with respect to the space variable and the time variable respectively. We introduce the initial data spaces $\widetilde{H}^\sigma(\mathbb{R})$, $\sigma\geq 0$:
\begin{align}\label{def4}
&\widetilde{H}^{\sigma}(\mathbb{R})=\Big\{\phi\in L^2(\mathbb{R}):\|\phi\|_{\widetilde{H}^{\sigma}}^2:=\|\eta_0\cdot\mathcal{F}_x(\phi)\|_{B_0}^2+\sum_{k=1}^\infty 2^{2\sigma k}\|\eta_k\cdot\mathcal{F}_x(\phi)\|_{L^2}^2<\infty\Big\},
\end{align}
where $\{\eta_k\}^\infty_{k=0}$ are the symbols of nonhomogeneous dyadic decomposition operators, and the Banach space $B_0(\mathbb{R})$ is defined by
\begin{align}\label{def4'}
B_0(\mathbb{R})=&\Big\{f\in L^2(\mathbb{R}):\, f \text{ supported in }[-2,2]\text{ and }\nonumber\\
&\|f\|_{B_0}:=\inf_{f=g+h}\|\mathcal{F}_\xi^{-1}(g)\|_{L^1_x}+\sum_{k'=-\infty}^{1}2^{-k'/2}\|\chi_{k'}\cdot h\|_{L^2_{\xi}}<\infty\Big\},
\end{align}
where $\{\chi_{k'}\}^{+\infty}_{k'=-\infty}$ are the symbols of homogeneous dyadic decomposition operators.
It is easy to see from the definitions that $\widetilde{H}^{\sigma}\hookrightarrow {H}^{\sigma}$, $\sigma\geq 0$. Moreover, from the scaling point of view, we have
\begin{equation}\label{scaling}
\|\phi_\lambda\|_{\widetilde{H}^\sigma}\leq C \|\phi\|_{\widetilde{H}^\sigma}\text{ for any }\lambda\in(0,1]\text{ and }\sigma\geq 0,
\end{equation}
where $\phi_\lambda(x):=\lambda\phi(\lambda x)$. In fact, the spaces $\widetilde{H}^{\sigma}$ are scaling critical for the low-frequency part, due to $\|\eta_0\cdot\mathcal{F}_x(\phi_\lambda)\|_{B_0}\sim \|\eta_0\cdot\mathcal{F}_x(\phi)\|_{B_0}$ for any $\lambda\in(0,1]$. Because of this, the inequality \eqref{scaling} could not be improved and we can only allow small initial data.

Let $\widetilde {H}^\infty(\mathbb{R})=\bigcap_{\sigma=0}^\infty \widetilde {H}^\sigma(\mathbb{R})$ with the induced metric. Let $S_\epsilon^\infty:\widetilde {H}^\infty(\mathbb{R})\to C([-1,1]:\widetilde {H}^\infty(\mathbb{R}))$ denote the nonlinear mapping that associates to any data $\phi\in \widetilde {H}^\infty$ the corresponding classical solution $u\in C([-1,1]:\widetilde {H}^\infty)$ of the initial value problem \eqref{BOB}. For any Banach space $V$ and $r>0$, let $B(r,V)$ denote the open
ball $\{v\in V:\|v\|_V<r\}$. Our main theorem states uniform local
well-posedness of
the BOB initial-value problem \eqref{BOB} for small data in $\widetilde{H}^{\sigma}$, $\sigma\geq 0$.

\begin{thm}\label{thm1}
(a) For any $\epsilon\in(0,1]$, there exists a constant $\delta>0$ with the property that for any $\phi\in B(\delta,\widetilde{H}^0)\cap \widetilde{H}^\infty$ there is a unique solution
\begin{equation*}
u^\epsilon=S_\epsilon^\infty(\phi)\in C([-1,1]:\widetilde{H}^\infty)
\end{equation*}
of the initial-value  problem \eqref{BOB}.

(b) For any $\phi\in B(\delta,\widetilde{H}^0)$, the mapping $\phi\to S_\epsilon^\infty(\phi)$ extends (uniquely) to a Lipschitz mapping
\begin{equation*}
S_\epsilon^0:B(\delta,\widetilde{H}^0)\to C([-1,1]:\widetilde{H}^0),
\end{equation*}
uniformly on $\epsilon\in(0,1]$ with the property that $S_\epsilon^0(\phi)$ is a solution of the initial-value problem \eqref{BOB}.

(c) For any $\sigma \in[0,\infty)$ we have the local Lipschitz bound which is independent of $\epsilon$
\begin{equation*}
\sup_{t\in[-1,1]}\|S_\epsilon^0(\phi)(t)-S_\epsilon^0(\phi')(t)\|_{\widetilde{H}^{\sigma}}\leq C(\sigma,R)\|\phi-\phi'\|_{\widetilde{H}^{\sigma}}
\end{equation*}
for any $R>0$ and $\phi,\phi'\in B(\delta,\widetilde{H}^0)\cap B(R,\widetilde{H}^{\sigma})$. As a consequence, the mapping $S^0$ restricts to a locally Lipschitz mapping
\begin{equation*}
S_\epsilon^{\sigma}:B(\delta,\widetilde{H}^0)\cap \widetilde{H}^{\sigma}\to C([-1,1]:\widetilde{H}^{\sigma}),
\end{equation*}
uniformly on $\epsilon\in(0,1]$.

(d) For any $\sigma \in[0,\infty)$, denote $\phi\to S^\sigma(\phi)$ the solution mapping of the initial-value problem \eqref{BO}, then we have the limit behavior
\begin{align*}
\lim_{\epsilon\rightarrow0}\|S^\sigma_\epsilon(\phi)-S^\sigma(\phi)\|_{C([-1,1];\widetilde{H}^\sigma)}=0.
\end{align*}
\end{thm}


{\bf Notations.} In the sequel $C$ will denote a universal positive constant which
can be different at each appearance. $x\lesssim y $ (for $x$, $y
>0$) means that $x\le Cy$, and $x\sim y$ stands for $x\lesssim y $
and $y\lesssim x$. $\mathscr{F}$ ($\mathscr{F}^{-1}$) denotes the (inverse) Fourier transform. $\widehat{\varphi}$ also denotes the Fourier transform of a distribution $\varphi$.

\section{Function spaces and known results }\label{section2}

At the beginning, let us recall the dyadic decomposition. Denote $\mathbb{Z}_+=\mathbb{Z}\cap[0,\infty)$.
Let $\eta_0:\mathbb{R}\to[0,1]$ denote an even smooth function
supported in $[-8/5,8/5]$ and equal to $1$ in $[-5/4,5/4]$. For
$\ell\in\mathbb{Z}$ let
$\chi_\ell(\xi)=\eta_0(\xi/2^\ell)-\eta_0(\xi/2^{\ell-1})$, $\chi_\ell$
supported in $\{\xi:|\xi|\in[(5/8)\cdot2^{\ell},(8/5)\cdot2^{\ell}]\}$, and
\begin{equation*}
\chi_{[\ell_1,\ell_2]}=\sum_{\ell=\ell_1}^{\ell_2}\chi_\ell\ \text{ for any }\ \ell_1\leq \ell_2\in\mathbb{Z}.
\end{equation*}
For simplicity of notation, let $\eta_\ell=\chi_\ell$ if $\ell\geq 1$ and
$\eta_\ell\equiv 0$ if $\ell\leq -1$. Also, for $\ell_1\leq
\ell_2\in\mathbb{Z}$ let
\begin{equation*}
\eta_{[\ell_1,\ell_2]}=\sum_{\ell=\ell_1}^{\ell_2}\eta_\ell\ \text{ and }\ \eta_{\leq \ell_2}=\sum_{\ell=-\infty}^{\ell_2}\eta_\ell.
\end{equation*}
For any $k\in \mathbb{Z}_+$ and $\phi\in L^2(\mathbb{R})$ we define
the operator $P_k$ by the formula
\begin{equation*}
\mathcal{F}_x(P_k\phi)(\xi)=\eta_k(\xi)\mathcal{F}_x(\phi)(\xi).
\end{equation*}
By a slight abuse of notation we also define the operators $P_k$ on $L^2(\mathbb{R}\times\mathbb{R})$ by the formula $\mathcal{F}(P_ku)(\xi,\tau)=\eta_k(\xi)\mathcal{F}(u)(\xi,\tau)$. For $\ell\in\mathbb{Z}$ let
$I_\ell=\{\xi\in\mathbb{R}:|\xi|\in[2^{\ell-1},2^{\ell+1}]\}$. For
$\ell\in \mathbb{Z}_+$ let $\widetilde{I}_\ell=[-2,2]$ if
$\ell=0$ and $\widetilde{I}_\ell=I_\ell$ if $\ell\geq 1$. For $k\in\mathbb{Z}$
and $j\geq 0$ let
\begin{equation*}
\begin{cases}
&D_{k,j}=\{(\xi,\tau)\in\mathbb{R}\times\mathbb{R}:\xi\in
I_k,\,\tau-\omega(\xi)\in\widetilde{I}_j\}\text{ if }k\geq 1;\\
&D_{k,j}=\{(\xi,\tau)\in\mathbb{R}\times\mathbb{R}:\xi\in
I_k,\,\tau\in\widetilde{I}_j\}\text{ if }k\leq 0.
\end{cases}
\end{equation*}
For $\xi\in\mathbb{R}$ let $\omega(\xi)$ denote the dispersive relation of BO equation, i.e.,
\begin{equation}\label{omega}
\omega(\xi)=-\xi|\xi|.
\end{equation}
\begin{defn} We define the Banach spaces $X_k=X_k(\mathbb{R}\times\mathbb{R})$,
$k\in\mathbb{Z}_+$: for $k\geq 1$ we define
\begin{equation}\label{def1}
\begin{split}
X_k=&\{f\in L^2:\, f \text{ supported in }I_k\times\mathbb{R}\text{ and }\\
&\|f\|_{X_k}:=\sum_{j=0}^\infty 2^{j/2}\beta_{k,j}\|\eta_j(\tau-\omega(\xi))f(\xi,\tau)\,\|_{L^2_{\xi,\tau}}<\infty\},
\end{split}
\end{equation}
where
\begin{equation}\label{def1'}
\beta_{k,j}=1+2^{(j-2k)/2}.
\end{equation}
For $k=0$ we define
\begin{equation}\label{def1''}
\begin{split}
X_0=&\{f\in L^2:\, f \text{ supported in }\widetilde{I}_0\times\mathbb{R}\text{ and }\\
&\|f\|_{X_0}:=\sum_{j=0}^\infty\sum_{k'=-\infty}^12^{j-k'/2}\|\eta_j(\tau)\chi_{k'}(\xi)f(\xi,\tau)\,\|_{L^2_{\xi,\tau}}<\infty\}.
\end{split}
\end{equation}
\end{defn}
The choices of the coefficients $\beta_{k,j}$ and the large factor $2^{-k'/2}$ are important in order to get the bilinear estimates.
The spaces $X_k$ are not sufficient for our purpose, due to
various logarithmic divergences involving the modulation variable. For $k\geq 100$ and $k=0$ we also define the
Banach spaces $Y_k=Y_k(\mathbb{R}\times\mathbb{R})$.
\begin{defn}
For $k\geq 100$ we define
\begin{equation}\label{def2}
\begin{split}
Y_k=\{f\in L^2:\,&f\text{ supported in } \bigcup_{j=0}^{k-1}D_{k,j}\text{ and }\\
&\|f\|_{Y_k}:=2^{-k/2}\|\mathcal{F}^{-1}[(\tau-\omega(\xi)+i)
f(\xi,\tau)]\|_{L^1_xL^2_t}<\infty\},
\end{split}
\end{equation}
where $i$ is the unit imaginary number.
For $k=0$ we define
\begin{equation}\label{def2''}
\begin{split}
Y_0=\{f\in L^2:\,&f\text{ supported in } \widetilde{I}_0\times\mathbb{R}\text{ and }\\
&\|f\|_{Y_0}:=\sum_{j\geq1}2^j\|\mathcal{F}^{-1}[\eta_j(\tau)
f(\xi,\tau)]\|_{L^1_xL^2_t}\\
&\quad\quad\quad\quad+\sum_{j\leq 0}\|\mathcal{F}^{-1}[\chi_{j}(\tau)
f(\xi,\tau)]\|_{L^1_xL^2_t}<\infty\}.
\end{split}
\end{equation}
\end{defn}
\begin{rem} The definition of $Y_0$ is different from that in \cite{Io-Ke,Io-Ke2}. It is easy to see that the space $Y_0$ in this paper is smaller than the corresponding space (denote it by $\bar{Y}_0$) in \cite{Io-Ke,Io-Ke2}, whose norm is given by
\begin{align*}
\|f\|_{\bar{Y}_0}:=\sum_{j=0}^\infty2^j\|\mathcal{F}^{-1}[\eta_j(\tau)
f(\xi,\tau)]\|_{L^1_xL^2_t}.
\end{align*}
We use the homogeneous dyadic decomposition to avoid the logarithmic divergences which occur in getting uniform estimates of Benjamin-Ono-Burgers equation.
\end{rem}
\begin{defn} We define our basic Banach spaces $Z_k$.
\begin{equation}\label{def3}
Z_k:=X_k\text{ if }1\leq k\leq 99\text{ and }Z_k:=X_k+Y_k\text{ if }k\geq 100\text{ or }k=0.
\end{equation}
\end{defn}

In some estimates we will also need the space $\overline{Z}_0$, $Z_0\subseteq\overline{Z}_0$.
\begin{defn}
\begin{equation}\label{def1'''}
\begin{split}
\overline{Z}_0=&\{f\in L^2(\mathbb{R}\times\mathbb{R}):\, f \text{ supported in }\widetilde{I}_0\times\mathbb{R}\text{ and }\\
&\|f\|_{\overline{Z}_0}:=\sum_{j=0}^\infty2^{j}\|\eta_j(\tau)f(\xi,\tau)\,\|_{L^2_{\xi,\tau}}<\infty\}.
\end{split}
\end{equation}
\end{defn}

For $k\in\mathbb{Z}_+$ let
\begin{equation*}
\begin{cases}
&A_k(\xi,\tau)=\tau-\omega(\xi)+i\ \ \ \ \text{ if }k\geq 1;\\
&A_k(\xi,\tau)=\tau+i\ \ \ \ \text{ if }k=0.\\
\end{cases}
\end{equation*}
\begin{defn}
For $\sigma\geq 0$ we define the Banach spaces
$F^{\sigma}=F^{\sigma}(\mathbb{R}\times\mathbb{R})$, and
$N^{\sigma}=N^{\sigma}(\mathbb{R}\times\mathbb{R})$:
\begin{equation}\label{def5}
\begin{split}
F^{\sigma}&=\Big\{u\in\mathcal{S}'(\mathbb{R}\times\mathbb{R}):\|u\|_{F^{\sigma}}^2:=\sum_{k=0}^\infty 2^{2\sigma k}
\|\eta_k(\xi)(I-\partial_\tau^2)\mathcal{F}(u)\|_{Z_k}^2
<\infty\Big\},
\end{split}
\end{equation}
and
\begin{equation}\label{def6}
\begin{split}
N^{\sigma}&=\Big\{u\in\mathcal{S}'(\mathbb{R}\times\mathbb{R}):\|u\|_{N^{\sigma}}^2:=\sum_{k=0}^\infty 2^{2\sigma k}\|\eta_k(\xi)
A_k(\xi,\tau)^{-1}\mathcal{F}(u)\|_{Z_k}^2<\infty\Big\}.
\end{split}
\end{equation}
\end{defn}

We establish some basic properties and known estimates which are similar to that in \cite{Io-Ke}. Using the definitions, if $k\geq 1$ and $f_k\in Z_k$ then $f_k$ can be written in the form
\begin{equation}\label{repr1}
\begin{cases}
&f_k=\sum\limits_{j=0}^\infty f_{k,j}+g_k;\\
&\sum\limits_{j=0}^\infty 2^{j/2}\beta_{k,j}\|f_{k,j}\|_{L^2}+\|g_k\|_{Y_k}\leq 2\|f_k\|_{Z_k},
\end{cases}
\end{equation}
such that $f_{k,j}$ is supported in $D_{k,j}$ and $g_k$ is supported in $\bigcup_{j=0}^{k-1} D_{k,j}$ (if $k\leq 99$ then $g_k\equiv 0$).
If $f_0\in Z_0$ then $f_0$ can be written in the form
\begin{equation}\label{repr2}
\begin{cases}
&f_0=\sum\limits_{j=0}^\infty \sum\limits_{k'=-\infty}^1f^{k'}_{0,j}+\sum\limits_{j=-\infty}^\infty g_{0,j};\\
&\sum\limits_{j=0}^\infty \sum\limits_{k'=-\infty}^12^{j-k'/2}\|f^{k'}_{0,j}\|_{L^2}+\sum\limits_{j\geq1} 2^j\|\mathcal{F}^{-1}(g_{0,j})\|_{L^1_xL^2_t}+\sum\limits_{j\leq0} \|\mathcal{F}^{-1}(g_{0,j})\|_{L^1_xL^2_t}\leq 2\|f_0\|_{Z_0},
\end{cases}
\end{equation}
such that $f^{k'}_{0,j}$ is supported in $D_{k',j}$ and $g_{0,j}$ is supported in $\widetilde{I}_0\times{I}_j$.

\begin{lem}\label{lem1}
(a) If $m,m':\mathbb{R}\to\mathbb{C}$, $k\geq 0$, and $f_k\in Z_k$ then
\begin{equation}\label{l1}
\begin{cases}
&\|m(\xi)f_k(\xi,\tau)\|_{Z_k}\leq C\|\mathcal{F}_\xi^{-1}(m)\|_{L^1(\mathbb{R})}\|f_k\|_{Z_k};\\
&\|m'(\tau)f_k(\xi,\tau)\|_{Z_k}\leq C\|m'\|_{L^\infty(\mathbb{R})}\|f_k\|_{Z_k}.
\end{cases}
\end{equation}

(b) If $k\geq 1$, $j\geq 0$, and $f_k\in Z_k$ then
\begin{equation}\label{l2}
\|\eta_j(\tau-\omega(\xi))f_k(\xi,\tau)\|_{X_k}\leq C\|f_k\|_{Z_k}.
\end{equation}

(c) If $k\geq 1$, $j\in[0,k]$, and $f_k$ is supported in $I_k\times\mathbb{R}$ then
\begin{equation}\label{l3}
\|\mathcal{F}^{-1}[\eta_{\leq j}(\tau-\omega(\xi))f_k(\xi,\tau)]\|_{L^1_xL^2_t}\leq C\|\mathcal{F}^{-1}(f_k)\|_{L^1_xL^2_t}.
\end{equation}
\end{lem}

\begin{lem}\label{lem2}
If $k\geq 0$, $t\in\mathbb{R}$, and $f_k\in Z_k$ then
\begin{equation}\label{l4}
\begin{cases}
&\big|\big|\int_\mathbb{R}f_k(\xi,\tau)e^{it\tau}\,d\tau\big|\big|_{L^2_\xi}
\leq C||f_k||_{Z_k}\text{ if }k\geq 1;\\
&\big|\big|\int_\mathbb{R}f_0(\xi,\tau)e^{it\tau}\,d\tau\big|\big|_{B_0}
\leq C||f_0||_{Z_0}\text{ if }k=0.
\end{cases}
\end{equation}
As a consequence,
\begin{equation}\label{l5}
F^\sigma\subseteq C(\mathbb{R}:\widetilde{H}^\sigma)\text{ for any }\sigma\geq 0.
\end{equation}
\end{lem}

\section{Uniform linear estimates}\label{linear}

In this section, we construct the uniform homogeneous and inhomogeneous linear estimates for BOB equation. The dissipative structure $-\epsilon \partial_x^2 u$ destroys some symmetries and brings some logarithmic divergences, which will bring several technical difficulties.

For $\phi\in L^2(\mathbb{R})$, let $W_\epsilon(t)\phi\in C(\mathbb{R};L^2)$ denote the solution of the free Benjamin-Ono-Burgers evolution given by
\begin{align}
W_\epsilon(t)\phi=\mathcal{F}_\xi^{-1}e^{it\omega(\xi)-t\epsilon\xi^2}\mathcal{F}_x\phi,
\end{align}
where $\omega(\xi)$ is defined in \eqref{omega}. Assume $\psi:\mathbb{R}\rightarrow[0,1]$ is an even smooth function
supported in $[-8/5,8/5]$ and equal to $1$ in $[-5/4,5/4]$. In the following discussions, the implicit constant in inequality sign ``\ $\lesssim$\ '' is independent of $\epsilon$. We first prove a uniform estimate for the free solution.
\begin{lem}\label{Lemma1}
If $\sigma\geq 0$ and $\phi\in \widetilde{H}^{\sigma}$ then for any $\epsilon\in [0,1]$,
\begin{equation*}
\|\psi(t)\cdot (W_\epsilon(t)\phi)\|_{F^{\sigma}}\leq C\|\phi\|_{\widetilde{H}^{\sigma}},
\end{equation*}
where the constant $C$ is independent of $\epsilon$.
\end{lem}
{\bf Proof.} It follows from the definition of $F^\sigma$ that
\begin{align*}
&\|\psi(t)\cdot (W_\epsilon(t)\phi)\|_{F^{\sigma}}^2=\sum_{k=0}^\infty 2^{2\sigma k}
\|\eta_k(\xi)(I-\partial_\tau^2)\mathcal{F}(\psi(t)W_\epsilon(t)\phi)\|_{Z_k}^2\\
\leq &\sum_{k\geq 1} 2^{2\sigma k}
\|\eta_k(\xi)(I-\partial_\tau^2)\mathcal{F}(\psi(t)W_\epsilon(t)\phi)\|_{X_k}^2+\|\eta_0(\xi)(I-\partial_\tau^2)\mathcal{F}(\psi(t) W_\epsilon(t)\phi)\|_{Z_0}^2.
\end{align*}
In view of the definition of $\widetilde{H}^{\sigma}$, it suffices to prove that
\begin{align}
\|\eta_0(\xi)(I-\partial_\tau^2)\mathcal{F}(\psi(t) W_\epsilon(t)\phi)\|_{Z_0}&\lesssim \|\eta_0(\xi)\mathcal{F}_x(\phi)(\xi)\|_{B_0};\label{k=0}\\
\|\eta_k(\xi)(I-\partial_\tau^2)\mathcal{F}(\psi(t)W_\epsilon(t)\phi)\|_{X_k}&\lesssim
\|\eta_k(\xi)\mathcal{F}_x(\phi)(\xi)\|_{L^2},\ \ \  \text{for\ any}\  k\geq 1.\label{kgeq1}
\end{align}
Denote $\varphi(t):=(1+t^2)\psi(t)\in \mathcal{S}(\mathbb{R}^+)$, we have
\begin{align}
  (I-\partial_\tau^2)\mathcal{F}(\psi(t)W_\epsilon(t)\phi)&=\mathcal{F}_t((1+t^2)\psi(t)e^{it\omega(\xi)-t\epsilon\xi^2})\mathcal{F}_x(\phi)(\xi)\nonumber\\
  &=(\mathcal{F}_t(\varphi(t)e^{-t\epsilon\xi^2}))(\tau-\omega(\xi))\mathcal{F}_x(\phi)(\xi).\label{denote}
\end{align}

(1) $k=0$, proof of \eqref{k=0}. From \eqref{denote} we have
\begin{align}
\|\eta_0(\xi)(I-\partial_\tau^2)\mathcal{F}(\psi(t) W_\epsilon(t)\phi)\|_{Z_0}=\|\eta_0(\xi)\mathcal{F}_x(\phi)(\xi)\mathcal{F}_t (\varphi(t)e^{-t\epsilon\xi^2})(\tau-\omega(\xi))\|_{Z_0}.\label{k=0-1}
\end{align}
Write $\eta_0\cdot\mathcal{F}_x(\phi)=g(\xi)+\sum_{k'\leq 1}h_{k'}$, where $h_{k'}$ is supported in $I_{k'}$, then
\begin{equation}
\|\mathcal{F}_{\xi}^{-1}(g)\|_{L^1_x}+\sum_{k'\leq 1}2^{-k'/2}\|h_{k'}\|_{L^2}\leq 2\|\eta_0\cdot\mathcal{F}_x(\phi)\|_{B_0},\label{k=0-3}
\end{equation}
and \eqref{k=0-1} is controlled by
\begin{align}
\|g(\xi)\mathcal{F}_t (\varphi(t)e^{-t\epsilon\xi^2})(\tau-\omega(\xi))\|_{Z_0}+\sum_{k'\leq1}\|h_{k'}(\xi)\mathcal{F}_t (\varphi(t)e^{-t\epsilon\xi^2})(\tau-\omega(\xi))\|_{X_0}. \label{k=0-4}
\end{align}

We divide the first term in \eqref{k=0-4} into two parts $I+II$ as follows
\begin{align*}
\|g(\xi)\mathcal{F}_t (\varphi(t)e^{-t\epsilon\xi^2})(\tau)\|_{Y_0}+\big\|g(\xi)\big[\mathcal{F}_t (\varphi(t)e^{-t\epsilon\xi^2})(\tau-\omega(\xi))-\mathcal{F}_t( \varphi(t)e^{-t\epsilon\xi^2})(\tau)\big]\big\|_{X_0}.
\end{align*}
For the term $I$, by the definition and Young's inequality, we know that
\begin{align}
&\|g(\xi)\mathcal{F}_t (\varphi(t)e^{-t\epsilon\xi^2})(\tau)\|_{Y_0}\nonumber\\
=&\sum_{j\geq1} 2^j\big\|\mathcal{F}_\xi^{-1}[g(\xi)\eta_j(\tau)\mathcal{F}_t (\varphi(t)e^{-t\epsilon\xi^2})(\tau)]\big\|_{L^1_xL^2_\tau}\nonumber\\
&+\sum_{j'\leq0} \big\|\mathcal{F}_\xi^{-1}[g(\xi)\chi_{j'}(\tau)\mathcal{F}_t (\varphi(t)e^{-t\epsilon\xi^2})(\tau)]\big\|_{L^1_xL^2_\tau}\nonumber\\
\lesssim& \sum_{j\geq1} 2^j\|\mathcal{F}_\xi^{-1}g(\xi)\|_{L^1_x}\big\|\mathcal{F}_\xi^{-1}[\eta_{[0,1]}(\xi)\eta_j(\tau)\mathcal{F}_t (\varphi(t)e^{-t\epsilon\xi^2})(\tau)]\big\|_{L^1_xL^2_\tau}\nonumber\\
&+\sum_{j'\leq0} \|\mathcal{F}_\xi^{-1}g(\xi)\|_{L^1_x}\big\|\mathcal{F}_\xi^{-1}[\eta_{[0,1]}(\xi)\chi_{j'}(\tau)\mathcal{F}_t( \varphi(t)e^{-t\epsilon\xi^2})(\tau)]\big\|_{L^1_xL^2_\tau}.\label{I1}
\end{align}
It suffices to prove that
\begin{align}\sum_{j\geq1} 2^j\big\|\mathcal{F}_\xi^{-1}[\eta_{[0,1]}(\xi)\eta_j(\tau)\mathcal{F}_t (\varphi(t)e^{-t\epsilon\xi^2})(\tau)]\big\|_{L^1_xL^2_\tau}\lesssim 1,\label{k=0-2}
\end{align}
and
\begin{align}\sum_{j'\leq0}\big\|\mathcal{F}_\xi^{-1}[\eta_{[0,1]}(\xi)\chi_{j'}(\tau)\mathcal{F}_t (\varphi(t)e^{-t\epsilon\xi^2})(\tau)]\big\|_{L^1_xL^2_\tau}\lesssim 1.\label{k=0-7}
\end{align}
We divide them into $|x|\leq C$ and $|x|> C$ two cases.
If $|x|\leq C$, by H\"older's inequality and Taylor's expansion we know that
\begin{align*}
&\sum_{j\geq1} 2^j\big\|\mathcal{F}_\xi^{-1}[\eta_{[0,1]}(\xi)\eta_j(\tau)\mathcal{F}_t (\varphi(t)e^{-t\epsilon\xi^2})(\tau)]\big\|_{L^1_{|x|\leq C}L^2_\tau}\nonumber\\
=& \sum_{j\geq1} 2^j\Big\|\eta_{[0,1]}(\xi)\eta_j(\tau)\mathcal{F}_t \Big(\varphi(t)\sum_{n\geq 0}\frac{(-1)^n\epsilon^n |\xi|^{2n}}{n!}t^n\Big)(\tau)\Big\|_{L^2_\xi L^2_\tau}\nonumber\\
\lesssim &\sum_{n\geq 0}\frac{C^n}{n!}\sum_{j\geq1} 2^j\|\eta_j(\tau)\mathcal{F}_t (\varphi(t)t^n)\|_{L^2_\tau}\lesssim \sum_{n\geq 0}\frac{C^n}{n!}\|\varphi(t)t^n\|_{H^2}\lesssim 1.
\end{align*}
Similarly, combining with Hausdorff-Young inequality, we can get
\begin{align*}
&\sum_{j'\leq0}\big\|\mathcal{F}_\xi^{-1}[\eta_{[0,1]}(\xi)\chi_{j'}(\tau)\mathcal{F}_t (\varphi(t)e^{-t\epsilon\xi^2})(\tau)]\big\|_{L^1_{|x|\leq C}L^2_\tau}\nonumber\\
=& \sum_{n\geq 0}\frac{C^n}{n!}\sum_{j'\leq0}\|\chi_{j'}(\tau)\|_{L^2_\tau}\|\mathcal{F}_t (\varphi(t)t^n)(\tau)\|_{L^\infty_\tau}\nonumber\\
\lesssim& \sum_{n\geq 0}\frac{C^n}{n!}\|\varphi(t)t^n\|_{L^1_t}\lesssim 1,
\end{align*}
where we used the fact that $\sum_{j'\leq0}\|\chi_{j'}(\tau)\|_{L^2_\tau}\lesssim \sum_{j'\leq0} 2^{j'/2}\lesssim 1$. If $|x|>C$, then $|x|\sim \langle x\rangle$.  For any fixed $x$, we have
\begin{align*}
&\sum_{j\geq1} 2^j\|\eta_j(\tau)\mathcal{F}_t [\varphi(t)\mathcal{F}^{-1}_{\xi}(e^{-t\epsilon\xi^2})](\tau)\|_{L^2_\tau}\lesssim \|\varphi(t)\mathcal{F}^{-1}_{\xi}(e^{-t\epsilon\xi^2})\|_{H^2_t}\\
&\quad\quad \lesssim\|\varphi(t)(\sqrt{t\epsilon})^{-1}e^{-|x|^2/(t\epsilon)}\|_{H^2_t}\lesssim |x|^{-2},
\end{align*}
and
\begin{align*}
&\sum_{j'\leq0}\|\chi_{j'}(\tau)\mathcal{F}_t [\varphi(t)\mathcal{F}^{-1}_{\xi}(e^{-t\epsilon\xi^2})](\tau)\|_{L^2_\tau}\lesssim \|\varphi(t)\mathcal{F}^{-1}_{\xi}e^{-t\epsilon\xi^2}\|_{L^1_t}\\
&\quad\quad \lesssim\|\varphi(t)(\sqrt{t\epsilon})^{-1}e^{-|x|^2/(t\epsilon)}\|_{L^1_t}\lesssim |x|^{-2}.
\end{align*}
Therefore, one can get the conclusion \eqref{k=0-2} and \eqref{k=0-7}. For the term $II$, by the definition, the mean value theorem, and Taylor's expansion, for some $\theta\in[0,1]$, we have
\begin{align}
&\big\|g(\xi)\big[\mathcal{F}_t (\varphi(t)e^{-t\epsilon\xi^2})(\tau-\omega(\xi))-\mathcal{F}_t (\varphi(t)e^{-t\epsilon\xi^2})(\tau)\big]\big\|_{X_0}\nonumber\\
\lesssim&\sum_{j\geq0}\sum_{k'\leq 1}2^{j-k'/2}\|\eta_j(\tau)\chi_{k'}(\xi)g(\xi)\xi^2\mathcal{F}_t (t\varphi(t)e^{-t\epsilon\xi^2})(\tau-\theta\omega(\xi))\,\|_{L^2_{\xi,\tau}}\nonumber\\
\lesssim& \sum_{j\geq0} 2^j \sup_{|\xi|\leq 2} \|P_j(t\varphi(t)e^{-t\epsilon\xi^2} e^{it\theta\omega(\xi)})\|_{L^2_t} \sum_{k'\leq 1} 2^{2k'}\|g(\xi)\|_{L^\infty_\xi}\nonumber\\
\lesssim& \sum_{j\geq0} 2^j  \|P_j(\varphi(t)t^{n+1})\|_{L^2_t}\|\mathcal{F}_{\xi}^{-1}(g)\|_{L^1_x} \sup_{|\xi|\leq 2} \sum_{n\geq 0}\frac{\big|\epsilon\xi^2+i\theta\xi|\xi|\big|^n}{n!}\nonumber\\
\lesssim& \sum_{n\geq 0}\frac{C^n}{n!}\|\varphi(t)t^{n+1}\|_{H^2}\|\mathcal{F}^{-1}g(\xi)\|_{L^1_x}\lesssim \|\mathcal{F}_{\xi}^{-1}(g)\|_{L^1_x}.\label{I2}
\end{align}
In view of \eqref{I1}-\eqref{I2}, we can get that
\begin{align}
\|g(\xi)\mathcal{F}_t (\varphi(t)e^{-t\epsilon\xi^2})(\tau-\omega(\xi))\|_{Z_0}\lesssim \|\mathcal{F}_{\xi}^{-1}(g)\|_{L^1_x}. \label{k=0-5}
\end{align}

For the second term in \eqref{k=0-4}, recall that  $h_{k'}$ is supported in $I_{k'}$, from the definition and Taylor's expansion,  we can obtain that for any fixed $k'$,
\begin{align}
&\|h_{k'}(\xi)\mathcal{F}_t (\varphi(t)e^{-t\epsilon\xi^2})(\tau-\omega(\xi))\|_{X_0}\nonumber\\
\lesssim&\sum_{j\geq0} 2^{j-k'/2}\|\eta_j(\tau)h_{k'}(\xi)\mathcal{F}_t (\varphi(t)e^{-t\epsilon\xi^2}e^{-it\xi|\xi|})(\tau)\,\|_{L^2_{\xi,\tau}}\nonumber\\
\lesssim& \sum_{j\geq0} 2^j \sup_{|\xi|\leq 2} \|P_j(\varphi(t)e^{-t\epsilon\xi^2} e^{-it\xi|\xi|})\|_{L^2_t} \cdot2^{-k'/2}\|h_{k'}(\xi)\|_{L^2_\xi}\nonumber\\
\lesssim& \sum_{j\geq0} 2^j \sup_{|\xi|\leq 2} \bigg\|P_j(\varphi(t)t^n) \sum_{n\geq 0}\frac{(-1)^n(\epsilon\xi^2+i\xi|\xi|)^n}{n!}\bigg\|_{L^2_t}\cdot2^{-k'/2}\|h_{k'}(\xi)\|_{L^2_\xi}\nonumber\\
\lesssim& \sum_{n\geq 0}\frac{C^n}{n!}\|\varphi(t)t^{n}\|_{H^2}\cdot2^{-k'/2}\|h_{k'}(\xi)\|_{L^2_\xi}\lesssim 2^{-k'/2}\|h_{k'}(\xi)\|_{L^2_\xi}. \label{k=0-6}
\end{align}
Therefore, combining \eqref{k=0-1}-\eqref{k=0-4} and \eqref{k=0-5}-\eqref{k=0-6}, we obtain the conclusion \eqref{k=0}.

(2) $k\geq 1$, proof of \eqref{kgeq1}. For any $k\geq 1$, by the change of variables and H\"older's inequality, we get
\begin{align*}
\|\eta_k(\xi)(I-\partial_\tau^2)\mathcal{F}(\psi(t)W_\epsilon(t)\phi)\|_{X_k}&=\sum_{j\geq0} 2^{j/2}\beta_{k,j}\|\eta_k(\xi)\mathcal{F}_x(\phi)(\xi)\eta_j(\tau)\mathcal{F}_t(\varphi(t)e^{-t\epsilon\xi^2})(\tau)\|_{L^2_{\xi,\tau}}\nonumber\\
&\lesssim \|\eta_k(\xi)\mathcal{F}_x(\phi)(\xi)\|_{L^2_{\xi}} \sum_{j\geq0} 2^{j/2}\beta_{k,j}\sup_{|\xi|\sim2^k}\|P_j(\varphi(t)e^{-t\epsilon\xi^2})\|_{L^2_{t}}.
\end{align*}
It suffices to show that for any $k\geq 1$,
\begin{equation}\label{claim}
\sum_{j\geq0} 2^{j/2}\beta_{k,j}\sup_{|\xi|\sim2^k}\|P_j(\varphi(t)e^{-t\epsilon\xi^2})\|_{L^2_{t}}\lesssim 1,
\end{equation}
where the implicit constant is independent of $\epsilon$ and $k$. By using Plancherel's equality and the fact that
$$\mathcal{F}_t(e^{-|t|})(\tau)=C\frac{1}{1+|\tau|^2},$$
we know that if $|\xi|\sim 2^k$, then for any $ j\geq 0$,
\begin{equation}\label{claim2}
\|P_j(e^{-\epsilon \xi^2|t|})(t)\|_{L^2}\lesssim
\|P_j(e^{-\epsilon 2^{2k}|t|})(t)\|_{L^2}.
\end{equation}
To prove \eqref{claim} we may assume $j\geq 100$ in the summation. Using the para-product homogeneous
decomposition, we have
\begin{align}
P_j(u_1u_2) = P_j \Big(\sum_{r\geq j-10} (P_{r+1} u_1) (P_{\leq
r}u_2)+ (P_{\leq r}u_1) (P_{r}u_2)\Big):=P_j(I+II).
\label{para}
\end{align}
Now we take $u_1=e^{-\epsilon
|t|\xi^2}$ and $u_2=\varphi(t)$.  For $P_j(I)$, it follows from
H\"older's inequality and \eqref{claim2} that
\begin{align*}
&\sum_{j \ge 100} 2^{j/2} \beta_{k,j}\sup_{|\xi|\sim2^k}\|P_j(I)\|_{L^2_t} \nonumber\\
\lesssim &
\sum_{j \ge 100} 2^{j/2} \beta_{k,j} \sum_{r\geq j-10} \sup_{|\xi|\sim2^k}\| P_{r+1}e^{-\epsilon
|t|\xi^2}\|_{L^{2}_t}
\|P_{\leq r}\varphi(t)\|_{L^\infty_{t}} \nonumber\\
\lesssim & \sum_{j\geq 100} \big(2^{j/2}+2^{j-k}\big) \sum_{r\geq j-10} \|
P_{r+1} e^{-\epsilon
|t|2^{2k}}\|_{L^{2}_t}
:=I_1+I_2.
\end{align*}
Then by discrete Young's inequality we can get
\begin{align*}
I_1\leq \sum_{j\geq 100} \sum_{r\geq j-10}2^{(j-r)/2}2^{r/2}\|
P_{r+1} e^{-\epsilon
|t|2^{2k}}\|_{L^{2}_t}\lesssim \|e^{-\epsilon
|t|2^{2k}}\|_{\dot{B}^{1/2}_{2,1}}\lesssim \|e^{-
|t|}\|_{\dot{B}^{1/2}_{2,1}}\lesssim 1,
\end{align*}
\begin{align*}
I_2\leq 2^{-k}\sum_{j\geq 100} \sum_{r\geq j-10}2^{j-r}2^{r}\|
P_{r+1} e^{-\epsilon
|t|2^{2k}}\|_{L^{2}_t}\lesssim 2^{-k}\|e^{-\epsilon
|t|2^{2k}}\|_{\dot{B}^{1}_{2,1}}\lesssim \epsilon^{1/2}\|e^{-
|t|}\|_{\dot{B}^{1}_{2,1}}\lesssim 1,
\end{align*}
where we used the facts that $\|e^{-\lambda
|t|}\|_{\dot{B}^{\sigma}_{2,1}}\sim \lambda^{\sigma-1/2}\|e^{-
|t|}\|_{\dot{B}^{\sigma}_{2,1}}$ and $e^{-|t|}\in\dot{B}_{2,1}^{1/2},\ \dot{B}_{2,1}^{1}$.
For $P_j(II)$, it follows from Bernstein's estimate, H\"older's inequality and \eqref{claim2} that
\begin{align*}
\sum_{j \ge 100} 2^{j/2} \beta_{k,j}\sup_{|\xi|\sim2^k}\|P_j(II)\|_{L^2_t}
\lesssim &
\sum_{j \ge 100} 2^{j/2} \beta_{k,j} \sum_{r\geq j-10} \|P_{r+1}\varphi(t)\|_{L^2_{t}}\sup_{|\xi|\sim2^k}\| P_{\leq r}e^{-\epsilon
|t|\xi^2}\|_{L^{\infty}_t}
 \nonumber\\
\lesssim & \|\varphi(t)\|_{\dot{B}^1_{2,1}}\sum_m 2^{m/2} \| P_{m}e^{-\epsilon
2^{2k}|t|}\|_{L^2_t}\lesssim \|e^{-
|t|}\|_{\dot{B}^{1/2}_{2,1}}\lesssim 1.
\end{align*}
Now we obtain the conclusion \eqref{claim} and then complete the proof of \eqref{kgeq1}. $\hfill\Box$

Before giving the inhomogeneous linear estimates, we state an important lemma, which will conquer the singularity when treating  $1/(\tau-\omega(\xi)-i\epsilon\xi^2)$. In addition, this lemma will effectively simplify the proof of uniform inhomogeneous estimates.
\begin{lem}\label{lemma3} If one of the following two assumptions holds:

(1) $k\geq 100$, $f_k$ is supported in $\bigcup_{j=0}^{k-1}D_{k,j}$ such that $f_k\in Y_k$;

(2) $k=0$, $f_0$ is supported in $\tilde{I}_0\times\mathbb{R}$ such that $f_0\in Y_0$,\\
then for any $\epsilon\in [0,1]$,
\begin{align*}
\frac{\tau-\omega(\xi)}{\tau-\omega(\xi)-i\epsilon\xi^2}f_k(\xi,\tau)\in Y_k, \ \ \ \ k\geq 100\ \text{or}\  k=0.
\end{align*}
In particular, we have
\begin{align}
\bigg\|\frac{\tau-\omega(\xi)}{\tau-\omega(\xi)-i\epsilon\xi^2}f_k(\xi,\tau)\bigg\|_{Y_k}\lesssim \|f_k\|_{Y_k}, \ \ \ \ k\geq 100\ \text{or}\ k=0.\label{L3-1}
\end{align}
\end{lem}
{\bf Proof.} (1) $k\geq100$. By the definition of $Y_k$, it suffices to prove that
\begin{align*}
\bigg\|\mathcal{F}^{-1}\frac{\tau-\omega(\xi)}{\tau-\omega(\xi)-i\epsilon\xi^2}(\tau-\omega(\xi)+i)f_k(\xi,\tau)\bigg\|_{L_x^1L_t^2}\lesssim \|\mathcal{F}^{-1}(\tau-\omega(\xi)+i)f_k\|_{L_x^1L_t^2}.
\end{align*}
In view of Plancherel's theorem and the support of $f_k$, we only need to prove that
\begin{align}
\bigg\|\int_{\mathbb{R}}e^{ix\xi}\frac{\tau-\omega(\xi)}{\tau-\omega(\xi)-i\epsilon\xi^2}\eta_{\leq k}(\tau-\omega(\xi))\chi_{[k-1,k+1]}(\xi)d\xi\bigg\|_{L_x^1L_\tau^\infty}\leq C. \label{L3}
\end{align}
The function in the left-hand side of \eqref{L3} is not zero only if $|\tau|\sim 2^{2k}$. By symmetry, we may assume $\{(\xi,\tau):\xi\in[2^{k-2}, 2^{k+2}],\tau\in[-2^{2k+10},-2^{2k-10}]\}$. Rewrite
\begin{align*}
\frac{\tau+\xi^2}{\tau+\xi^2-i\epsilon\xi^2}=\frac{1}{1-i\epsilon}\bigg(1+\frac{-i\epsilon\tau}{\tau+\xi^2-i\epsilon\xi^2}\bigg):=I+II.
\end{align*}

For $I$, by integration by parts, it is easy to show that
\begin{align*}
\bigg|\int\frac{(I-\Delta_\xi)(e^{ix\xi})}{1+x^2}\eta_{\leq k}(\tau+\xi^2)\chi_{[k-1,k+1]}(\xi)d\xi\bigg|\lesssim \frac{1}{1+x^2},
\end{align*}
where we used the fact $|\{\xi\in[2^{k-2}, 2^{k+2}]: |\tau+\xi^2|\leq 2^{k+2}\}|\leq C$.

For $II$, the case $|x|\leq 1$ is trivial, thus we just consider $|x|\geq 1$. Indeed, let
$$a:=\sqrt{\frac{\tau}{1-i\epsilon}}=-i|\tau|^{1/2}(1+\epsilon^2)^{-1/4}\Big(\cos \frac{\arctan\epsilon}{2} + i \sin \frac{\arctan\epsilon}{2}\Big), \ \ Re\ a\sim |\tau|^{1/2}\epsilon,$$
then by the fact
\begin{align*}
a\cdot\mathcal{F}^{-1}_\xi\frac{1}{\xi^2+a^2}=Ce^{-a|x|}, \ \ \ Re\ a>0,
\end{align*}
we can get
\begin{align*}
\mathcal{F}^{-1}_\xi\bigg(\frac{-i\epsilon\tau}{\tau+\xi^2-i\epsilon\xi^2}\bigg)=\frac{-i\epsilon\tau}{1-i\epsilon}\mathcal{F}^{-1}_\xi
\bigg(\frac{1}{\frac{\tau}{1-i\epsilon}+\xi^2}\bigg)=-i\epsilon C a\cdot
e^{-a|x|},
\end{align*}
Therefore,
\begin{align*}
\bigg|\mathcal{F}^{-1}_\xi\bigg(\frac{-i\epsilon\tau}{\tau+\xi^2-i\epsilon\xi^2}\bigg)\bigg|\leq C \epsilon^2|\tau|^{1/2}e^{-c\epsilon|\tau|^{1/2}|x|}\lesssim \epsilon2^ke^{-c\epsilon2^k|x|},
\end{align*}
whose $L^1_x$ norm is bounded, then we get the conclusion \eqref{L3}.

(2) $k=0$. By the definition of $Y_0$, we need to show that for any $j\in\mathbb{Z}$,
\begin{align*}
\bigg\|\mathcal{F}^{-1}\frac{\tau-\omega(\xi)}{\tau-\omega(\xi)-i\epsilon\xi^2}\chi_j(\tau)f_0(\xi,\tau)\bigg\|_{L_x^1L_t^2}\lesssim \|\mathcal{F}^{-1}\chi_j(\tau)f_0(\xi,\tau)\|_{L_x^1L_t^2}.
\end{align*}
Combining the Plancherel's theorem with Young's inequality, it suffices to prove that
\begin{align}\label{k=00}
\bigg\|\mathcal{F}_{\xi}^{-1}\frac{\tau-\omega(\xi)}{\tau-\omega(\xi)-i\epsilon\xi^2}\chi_j(\tau)\eta_{[0,1]}(\xi)\bigg\|_{L_x^1L_\tau^\infty}\lesssim 1.
\end{align}
Similar to (1), we may assume $\xi>0$ and rewrite
\begin{align*}
\frac{\tau+\xi^2}{\tau+\xi^2-i\epsilon\xi^2}=\frac{1}{1-i\epsilon}\bigg(1+\frac{-i\epsilon\tau}{\tau+\xi^2-i\epsilon\xi^2}\bigg):=I+II.
\end{align*}
Notice that
\begin{align*}
\bigg|\mathcal{F}^{-1}_\xi\bigg(\frac{-i\epsilon\tau}{\tau+\xi^2-i\epsilon\xi^2}\bigg)\bigg|\leq C \epsilon^2|\tau|^{1/2}e^{-c\epsilon|\tau|^{1/2}|x|}\lesssim \epsilon2^{j/2}e^{-c\epsilon2^{j/2}|x|}\in L^1_x,
\end{align*}
then we can get \eqref{k=00} in the same way as we used in (1). The proof is completed. $\hfill\Box$

For the inhomogeneous linear operator, we have the following uniform estimates.
\begin{lem}\label{Lemma2}
If $\sigma\geq 0$ and $u\in N^{\sigma}$,
then for any $\epsilon\in [0,1]$,
\begin{equation*}
\Big\|\psi(t)\cdot \int_0^tW_\epsilon(t-s)(u(s))\,ds\Big\|_{F^{\sigma}}\leq C\|u\|_{N^{\sigma}},
\end{equation*}
where the constant $C$ is independent of $\epsilon$.
\end{lem}
{\bf Proof.} By the definitions, it suffices to prove that $\forall\ k\geq 0$,
\begin{equation}\label{L2-1}
  \Big\|\eta_k(\xi)(I-\partial_\tau^2)\mathcal{F}\Big[\psi(t)\cdot \int_0^tW_\epsilon(t-s)(u(s))\,ds\Big]\Big\|_{Z_k}\lesssim \|\eta_k(\xi)
A_k(\xi,\tau)^{-1}\mathcal{F}(u)\|_{Z_k}.
\end{equation}
From a straightforward calculation, we have
\begin{align}\label{L2-2}
&(I-\partial_\tau^2)\mathcal{F}\Big[\psi(t)\cdot \int_0^tW_\epsilon(t-s)(u(s))ds\Big](\xi,\tau)\nonumber\\
=&\mathcal{F}_t\Big[(1+t^2)\psi(t)\cdot \int_0^t e^{-(t-s)\epsilon\xi^2}e^{-is \omega(\xi)}\mathcal{F}_x (u) (\xi,s) ds\Big](\tau-\omega(\xi))\nonumber\\
=&\mathcal{F}_t\Big[(1+t^2)\psi(t)e^{-t\epsilon\xi^2}\cdot \int_0^t e^{s\epsilon\xi^2}e^{-is \omega(\xi)}\int_\mathbb{R}e^{is\tau'}\mathcal {F}(u)(\xi,\tau')d\tau' ds\Big](\tau-\omega(\xi))\nonumber\\
=&\mathcal{F}_t\Big[(1+t^2)\psi(t)\cdot \int_\mathbb{R}\frac{e^{it(\tau'-\omega(\xi))}-e^{-t\epsilon\xi^2}}{i(\tau'-\omega(\xi))+\epsilon\xi^2}\mathcal {F} (u)(\xi,\tau')d\tau' \Big](\tau-\omega(\xi)).
\end{align}
Let $\varphi(t):=(1+t^2)\psi(t)$, and $f_k(\xi,\tau):=\eta_k(\xi)A_k(\xi,\tau)^{-1}\mathcal{F}(u)(\xi,\tau)$ for $k\in\mathbb{Z}_+$.
For $f_k\in Z_k$ let
\begin{equation}\label{L2-3}
T(f_k)(\xi,\tau):=\mathcal{F}_t\Big[\varphi(t)\cdot \int_\mathbb{R}\frac{e^{it(\tau'-\omega(\xi))}-e^{-t\epsilon\xi^2}}{i(\tau'-\omega(\xi))+\epsilon\xi^2}
A_k(\xi,\tau')f_k(\xi,\tau')d\tau' \Big](\tau-\omega(\xi)).
\end{equation}
In view of \eqref{L2-2}-\eqref{L2-3}, to prove \eqref{L2-1}, we only need to prove that
\begin{equation}\label{L2-4}
\|T\|_{Z_k\to Z_k}\leq C\ \ \text{ uniformly in }k\in\mathbb{Z}_+ \text{ and } \epsilon\in [0,1].
\end{equation}

(1) Case $k\geq 1$.

(1-$a$) Assume first that $f_k\in X_k$. The idea of this part is essential due to \cite{MR} and \cite{Guo-Wang}. Denote $f_k^\#(\xi,\mu')=f_k(\xi,\mu'+\omega(\xi))$ and $T(f_k)^\#(\xi,\mu)=T(f_k)(\xi,\mu+\omega(\xi))$. Then,
\begin{equation}\label{L2-5}
T(f_k)^\#(\xi,\mu)=\mathcal{F}_t\Big[\varphi(t)\cdot \int_\mathbb{R}\frac{e^{it\mu'}-e^{-t\epsilon\xi^2}}{i\mu'+\epsilon\xi^2}
(\mu'+i)f_k^\#(\xi,\mu')d\mu' \Big](\mu).
\end{equation}
It suffices to prove that
\begin{equation}\label{L2-6}
  \sum_{j\geq0} 2^{j/2}\beta_{k,j}\|\eta_j(\mu)T(f_k)^\#(\xi,\mu)\|_{L^2_{\xi,\mu}}\lesssim \sum_{j\geq0} 2^{j/2}\beta_{k,j}\|\eta_j(\mu)f_k^\#(\xi,\mu)\|_{L^2_{\xi,\mu}}.
\end{equation}
We divide $T(f_k)^\#(\xi,\mu)$ into four parts:
\begin{align*}
T(f_k)^\#(\xi,\mu)=&\mathcal{F}_t\Big[\varphi(t)\int_{|\mu'|\leq
1}\frac{e^{it\mu'}-1}{i\mu'+\epsilon\xi^2}(\mu'+i)f_k^\#(\xi,\mu')d\mu' \Big](\mu)\nonumber\\
&+\mathcal{F}_t\Big[\varphi(t)\int_{|\mu'|\leq
1}\frac{1-e^{-\epsilon
t\xi^2}}{i\mu'+\epsilon\xi^2}(\mu'+i)f_k^\#(\xi,\mu')d\mu' \Big](\mu)\nonumber\\
&+\mathcal{F}_t\Big[\varphi(t)\int_{|\mu'|\geq
1}\frac{e^{it\mu'}}{i\mu'+\epsilon\xi^2}(\mu'+i)f_k^\#(\xi,\mu')d\mu' \Big](\mu)\nonumber\\
&-\mathcal{F}_t\Big[\varphi(t)\int_{|\mu'|\geq
1}\frac{e^{-\epsilon
t\xi^2}}{i\mu'+\epsilon\xi^2}(\mu'+i)f_k^\#(\xi,\mu')d\mu' \Big](\mu)\nonumber\\
:=&I+II+III-IV.
\end{align*}
When $|\mu'|\geq 1$, the denominator in the fraction is far from 0, then $(\mu'+i)/(i\mu'+\epsilon\xi^2)$ is bounded, see the parts $III$ and $IV$. When $|\mu'|\leq 1$, we could use Taylor's expansion for the numerator to cancel the denominator, see the parts $I$ and $II$.
We now estimate the contributions of $I-IV$. Firstly, we consider the
contribution of $IV$.
\begin{align*}
&\sum_{j\geq0}2^{j/2}\beta_{k,j}\Big\|\eta_j(\mu)\mathcal{F}_t(\varphi(t)e^{-\epsilon
t\xi^2})(\mu)\int_{|\mu'|\geq
1}\frac{\mu'+i}{i\mu'+\epsilon\xi^2}f_k^\#(\xi,\mu')d\mu'\Big\|_{L_{\xi,\mu}^2}\nonumber\\
\lesssim&
\sum_{j\geq0}2^{j/2}\beta_{k,j}\sup_{\xi\in I_k}\|P_j(\varphi(t)e^{-\epsilon t\xi^2})(t)\|_{L_{t}^2}\cdot\int_{|\mu'|\geq 1}\|f_k^\#(\xi,\mu')\|_{L_\xi^2}d\mu'\nonumber\\
\lesssim&\int_{|\mu'|\geq 1}\|f_k^\#(\xi,\mu')\|_{L_\xi^2}d\mu'\lesssim \sum_{j\geq0}2^{j/2}\|\eta_j(\mu)f_k^\#(\xi,\mu)\|_{L^2_{\xi,\mu}},
\end{align*}
where we use the inequality \eqref{claim}. Secondly, we consider the contribution of $III$.
\begin{align*}
&\sum_{j\geq0}2^{j/2}\beta_{k,j}\Big\|\eta_j(\mu)\mathcal{F}_t\Big[\varphi(t)\int_{|\mu'|\geq
1}e^{it\mu'}\frac{\mu'+i}{i\mu'+\epsilon\xi^2}f_k^\#(\xi,\mu')d\mu'\Big](\mu)\Big\|_{L_{\xi,\mu}^2}\nonumber\\
\lesssim&\sum_{j\geq 0} 2^{j/2}\beta_{k,j}\Big\|\eta_j(\mu)|\mathcal{F}_t\varphi|\star\|f_k^\#(\xi,\cdot)\|_{L^2_\xi}\Big\|_{L^2_\mu}\nonumber\\
\lesssim&\sum_{j\geq 0} 2^{j/2}\beta_{k,j}\|P_j[\mathcal{F}^{-1}_\mu|\mathcal{F}_t\varphi|\cdot\mathcal{F}^{-1}_\mu\|f_k^\#(\xi,\mu)\|_{L^2_\xi}]
\|_{L^2_t}\nonumber\\
\lesssim& \sum_{j\geq0} 2^{j/2}\beta_{k,j}\|\eta_j(\mu)f_k^\#(\xi,\mu)\|_{L^2_{\xi,\mu}},
\end{align*}
where we used the facts that $B_{2,1}^{1/2}$ and $B_{2,1}^{1}$ are multiplication
algebras and that $\mathcal{F}_\mu^{-1}(|\mathcal{F}_t\varphi|)\in B_{2,1}^{1/2}$ and $B_{2,1}^{1}$.
Thirdly, we consider the contribution of $I$. By Taylor's expansion, we obtain
\begin{align*}
&\sum_{j\geq0}2^{j/2}\beta_{k,j}\Big\|\eta_j(\mu)\mathcal{F}_t\Big[\varphi(t)\int_{|\mu'|\leq
1}\sum_{n\geq
1}\frac{(it\mu')^n}{n!(i\mu'+\epsilon\xi^2)}(\mu'+i)f_k^\#(\xi,\mu')d\mu'\Big](\mu)\Big\|_{L_{\xi,\mu}^2}\nonumber\\
\lesssim&
\sum_{n\geq 1}\frac{\|t^n\varphi(t)\|_{B^{1}_{2,1}}}{n!}
\Big\|\int_{|\mu'|\leq1}\frac{|\mu'|}{|i\mu'+\epsilon\xi^2|}|f_k^\#(\xi,\mu')|d\mu'\Big\|_{L_\xi^2}\nonumber\\
\lesssim& \sum_{j\geq0} 2^{j/2}\|\eta_j(\mu)f_k^\#(\xi,\mu)\|_{L^2_{\xi,\mu}},
\end{align*}
where in the last inequality we used the fact
$\|t^n\varphi(t)\|_{B^{1}_{2,1}}\leq\|t^n\varphi(t)\|_{H^{2}}\leq C^n$.  Finally, we consider the contribution of $II$.
For $\epsilon\xi^2\geq 1$, the denominator in the fraction is far from 0, we can easily get that
\begin{align*}
&\sum_{j\geq0}2^{j/2}\beta_{k,j}\|\eta_j(\mu)II\|_{L_{\xi,\mu}^2}\nonumber\\
\lesssim&\sum_{j\geq0}2^{j/2}\beta_{k,j}\sup_{\xi\in I_k}\|P_j(\varphi(t)(1-e^{-\epsilon t\xi^2}))\|_{L_{t}^2}\cdot\int_{|\mu'|\leq 1}\|f_k^\#(\xi,\mu')\|_{L_\xi^2}d\mu'\nonumber\\
\lesssim& \sum_{j\geq0}2^{j/2}\|\eta_j(\mu)f_k^\#(\xi,\mu)\|_{L^2_{\xi,\mu}},
\end{align*}
where we use the inequality \eqref{claim} and $\varphi\in B^1_{2,1}$.
For $\epsilon\xi^2\leq 1$, using Taylor's expansion, we have
\begin{align*}
&\sum_{j\geq0}2^{j/2}\beta_{k,j}\|\eta_j(\mu)II\|_{L_{\xi,\mu}^2}\nonumber\\
\lesssim&\sum_{j\geq0}2^{j/2}\beta_{k,j}\Big\|\eta_j(\mu)\mathcal{F}_t\Big[\varphi(t)\sum_{n\geq
1}\frac{ t^n(\epsilon\xi^2)^n}{n!}\int_{|\mu'|\leq
1}\frac{(\mu'+i)}{(i\mu'+\epsilon\xi^2)}f_k^\#(\xi,\mu')d\mu'\Big]\Big\|_{L_{\xi,\mu}^2}\nonumber\\
\lesssim&
\sum_{n\geq 1}\frac{\|t^n\varphi(t)\|_{B^{1}_{2,1}}}{n!}
\Big\|\int_{|\mu'|\leq1}\frac{\epsilon\xi^2}{|i\mu'+\epsilon\xi^2|}|f_k^\#(\xi,\mu')|d\mu'\Big\|_{L_\xi^2}\nonumber\\
\lesssim& \sum_{j\geq0} 2^{j/2}\|\eta_j(\mu)f_k^\#(\xi,\mu)\|_{L^2_{\xi,\mu}}.
\end{align*}
Now we have shown that
\begin{equation}\label{L2-7}
\|T\|_{X_k\to X_k}\leq C\ \ \ \ \text{ uniformly in }k\geq 1 \ \text{and}\ \epsilon\in [0,1].
\end{equation}

(1-$b$) Assume now that $k\geq 100$, $f_k=g_k\in Y_k$. From \eqref{l2}, we know that  $\|\eta_j(\tau-\omega(\xi))f_k\|_{X_k}\lesssim \|f_k\|_{Y_k}$, thus we may assume that $g_k$
is supported in the set $\{(\xi,\tau):|\tau-\omega(\xi)|\leq
2^{k-20}\}$. For convenience, we decompose
\begin{equation}\label{gk}
g_k(\xi,\tau')=\frac{\tau'-\omega(\xi)}{\tau'-\omega(\xi)+i}g_k(\xi,\tau')+\frac{i}{\tau'-\omega(\xi)+i}g_k(\xi,\tau'),
\end{equation}
then \eqref{L2-3} becomes
\begin{align}\label{L2-11}
T(g_k)(\xi,\tau)=&\mathcal{F}_t\Big[\varphi(t)\cdot \int_\mathbb{R}\frac{e^{it(\tau'-\omega(\xi))}-e^{-t\epsilon\xi^2}}{i(\tau'-\omega(\xi))+\epsilon\xi^2}
(\tau'-\omega(\xi))g_k(\xi,\tau')d\tau' \Big](\tau-\omega(\xi))\nonumber\\
&+i T\Big( \frac{1}{\tau'-\omega(\xi)+i}g_k\Big)(\xi,\tau)\nonumber\\
=& \int_\mathbb{R}\frac{\tau'-\omega(\xi)}{i(\tau'-\omega(\xi))+\epsilon\xi^2}g_k(\xi,\tau')\hat{\varphi}(\tau-\tau')\,d\tau'\nonumber\\
&- \mathcal{F}_t(\varphi(t)e^{-t\epsilon\xi^2})(\tau-\omega(\xi))\int_\mathbb{R}\frac{\tau'-\omega(\xi)}{i(\tau'-\omega(\xi))+\epsilon\xi^2}g_k(\xi,\tau')\,d\tau'\nonumber\\
&+i T\Big( \frac{1}{\tau'-\omega(\xi)+i}g_k\Big)(\xi,\tau).
\end{align}

We can use \eqref{L2-7} to control the third term in \eqref{L2-11}. Notice that
\begin{align}\label{measure}
|\{\xi\in
I_k:|\tau-\omega(\xi)|\leq 2^{j+1}\}|\lesssim 2^{j-k},
\end{align}
then we have from \eqref{L2-7} that
\begin{align*}
 &\|T\big((\tau'-\omega(\xi)+i)^{-1}g_k\big)\|_{X_k}\lesssim\|(\tau'-\omega(\xi)+i)^{-1}g_k\|_{X_k}\nonumber\\
  \lesssim&\sum_{0\leq j\leq k}2^{j/2}\beta_{k,j}\|\eta_j(\tau'-\omega(\xi))(\tau'-\omega(\xi)+i)^{-1}g_k(\xi,\tau')\|_{L_{\xi,\tau'}^2}\nonumber\\
  \lesssim&\sum_{0\leq j\leq k}2^{-3j/2}2^{(j-k)/2}\|(\tau'-\omega(\xi)+i)g_k(\xi,\tau')\|_{L^\infty_\xi L^2_{\tau'}}
  \lesssim\|g_k\|_{Y_k}.
\end{align*}

For the first and second terms in \eqref{L2-11}, it suffices to prove that
\begin{align}\label{L2-a}
\Big\|\int_\mathbb{R}\frac{\tau'-\omega(\xi)}{\tau'-\omega(\xi)-i\epsilon\xi^2}g_k(\xi,\tau')\hat{\varphi}(\tau-\tau')\,d\tau'\Big\|_{Z_k}\lesssim \|g_k\|_{Y_k},
\end{align}
and
\begin{align}\label{L2-b}
\Big\|\mathcal{F}_t(\varphi(t)e^{-t\epsilon\xi^2})(\tau-\omega(\xi))\int_\mathbb{R}\frac{\tau'-\omega(\xi)}{\tau'-\omega(\xi)-i\epsilon\xi^2}g_k(\xi,\tau')\,d\tau'
\Big\|_{X_k}\lesssim \|g_k\|_{Y_k}.
\end{align}
Thanks to Lemma \ref{lemma3}, we know that
\begin{align*}
\Big\|\frac{\tau-\omega(\xi)}{\tau-\omega(\xi)-i\epsilon\xi^2}f_k(\xi,\tau)\Big\|_{Y_k}\lesssim \|f_k\|_{Y_k},
\end{align*}
then we can make the proof clearer and simpler. To prove \eqref{L2-a} and \eqref{L2-b}, we just need to prove
\begin{align}\label{L2-9}
\Big\|\int_\mathbb{R}g_k(\xi,\tau')\hat{\varphi}(\tau-\tau')\,d\tau'\Big\|_{Z_k}\lesssim \|g_k\|_{Y_k},
\end{align}
and
\begin{align}\label{L2-10}
\Big\|\mathcal{F}_t(\varphi(t)e^{-t\epsilon\xi^2})(\tau-\omega(\xi))\int_\mathbb{R}g_k(\xi,\tau')\,d\tau'
\Big\|_{X_k}\lesssim \|g_k\|_{Y_k}.
\end{align}

The inequality \eqref{L2-9} has been obtained by Ionescu and Kenig in \cite{Io-Ke}. For the sake of completeness, we give the rigorous proof. For the low modulation part, we divide it into two subparts:
\begin{equation*}
g_k(\xi,\tau')=g_k(\xi,\tau')\Big[\frac{\tau'-\omega(\xi)+i}{\tau-\omega(\xi)+i}+\frac{\tau-\tau'}{\tau-\omega(\xi)+i}\Big].
\end{equation*}
Then the left-hand side of \eqref{L2-9} is dominated by
\begin{align*}
&\Big\|\eta_{[0,k-1]}(\tau-\omega(\xi))(\tau-\omega(\xi)+i)^{-1}\int_\mathbb{R}g_k(\xi,\tau')(\tau'-\omega(\xi)+i)\hat{\varphi}(\tau-\tau')\,d\tau'\Big\|_{Y_k}\nonumber\\
&+\sum_{j\leq k-1}2^{j/2}\Big\|\eta_{j}(\tau-\omega(\xi))(\tau-\omega(\xi)+i)^{-1}\int_\mathbb{R}g_k(\xi,\tau')\hat{\varphi}(\tau-\tau')(\tau-\tau')\,d\tau'\Big\|_{L^2_{\xi,\tau}}\nonumber\\
&+\sum_{j\geq k-1}2^{j/2}\beta_{k,j}\Big\|\eta_{j}(\tau-\omega(\xi))\int_\mathbb{R}g_k(\xi,\tau')\hat{\varphi}(\tau-\tau')\,d\tau'\Big\|_{L^2_{\xi,\tau}}:=I+II+III.
\end{align*}
For $I$, we use Lemma \ref{lem1} (c) to bound it by
\begin{align*}
I&\lesssim 2^{-k/2}\|\mathcal{F}^{-1}\eta_{[0,k-1]}(\tau-\omega(\xi))\int_\mathbb{R}g_k(\xi,\tau')
(\tau'-\omega(\xi)+i)\hat{\varphi}(\tau-\tau')\,d\tau'\|_{L^1_xL^2_t}\nonumber\\ &\lesssim2^{-k/2}\|\varphi(t)\cdot\mathcal{F}^{-1}[(\tau'-\omega(\xi)+i)g_k(\xi,\tau')]\|_{L^1_xL^2_t}\lesssim\|g_k\|_{Y_k},
\end{align*}
as desired. For $II$, from \eqref{measure} we can get
\begin{align*}
  II\lesssim&\sum_{0\leq j\leq k}2^{-j/2}\Big\|\int_\mathbb{R}g_k(\xi,\tau')\hat{\varphi}(\tau-\tau')(\tau-\tau')\,d\tau'\Big\|_{L^2_{\xi,\tau}}\nonumber\\
  \lesssim&\|g_k(\xi,\tau)\|_{L^2_{\xi,\tau}}\lesssim\sum_{0\leq j\leq k-1}2^{-j}\|\eta_j(\tau-\omega(\xi))(\tau-\omega(\xi)+i)g_k(\xi,\tau)\|_{L_{\xi,\tau}^2}\nonumber\\
  \lesssim&\sum_{0\leq j\leq k-1}2^{-j}2^{(j-k)/2}\|(\tau-\omega(\xi)+i)g_k(\xi,\tau)\|_{L^\infty_\xi L^2_{\tau}}
  \lesssim\|g_k\|_{Y_k}.
\end{align*}
For $III$, let $g_k^\#(\xi,\mu')=g_k(\xi,\mu'+\omega(\xi))$, then
\begin{align*}
  III\lesssim&\sum_{j\geq k-1}2^{j}\Big\|\eta_{j}(\mu)\int_\mathbb{R}g_k^\#(\xi,\mu')\hat{\varphi}(\mu-\mu')\,d\mu'\Big\|_{L^2_{\xi,\mu}}\nonumber\\
  \lesssim&\sum_{j\geq k-1}2^{j}\sum_{j'\leq k-20}\Big\|\eta_{j}(\mu)\langle\mu\rangle^{-2}\int_\mathbb{R}\eta_{j'}(\mu')g_k^\#(\xi,\mu')\hat{\varphi}(\mu-\mu')\langle\mu-\mu'\rangle^2\,d\mu'\Big\|_{L^2_{\xi,\mu}}\nonumber\\
  \lesssim&\sum_{j'\leq k-20}\|\eta_{j'}(\mu)g_k^\#(\xi,\mu)\|_{L^2_{\xi,\mu}}\nonumber\\
  \lesssim&\sum_{j\leq k-20}2^{-j}\|\eta_j(\tau-\omega(\xi))(\tau-\omega(\xi)+i)g_k(\xi,\tau)\|_{L_{\xi,\tau}^2}
  \lesssim\|g_k\|_{Y_k}.
\end{align*}

Finally, to prove \eqref{L2-10}, we define the modified Hilbert
transform operator
\begin{equation*}
\mathcal{L}_k(g)(\mu):=\int_\mathbb{R}g(\tau)(\tau-\mu+i)^{-1}
\eta_{[0,k]}(\tau-\mu)\,d\tau,\,\,g\in L^2(\mathbb{R}).
\end{equation*}
Notice that
$$\mathcal{F}_{\tau}^{-1}\Big(\frac{1}{\tau-i}\Big)=\mathcal{F}_{\tau}^{-1}\Big(\frac{\tau}{\tau^2+1}\Big)+i\mathcal{F}_{\tau}^{-1}\Big(\frac{1}{\tau^2+1}\Big)=C(-isgn(t)+i)e^{-c t} \in L^\infty(\mathbb{R}),$$
Hence by Plancherel's theorem and H\"older's inequality, we have $\|\mathcal{L}_k\|_{L^2\to L^2}\leq C$, uniformly in $k$. We notice that if $g_k\in Y_k$ then $g_k$ can be
written in the form
\begin{equation*}
\begin{cases}
&g_k(\xi,\tau)=2^{k/2}\chi_{[k-1,k+1]}(\xi)(\tau-\omega(\xi)+i)^{-1}
\eta_{\leq k}(\tau-\omega(\xi))\int_\mathbb{R}e^{-ix\xi}h(x,\tau)\,dx;\\
&\|g_k\|_{Y_k}= C\|h\|_{L^1_xL^2_\tau}.
\end{cases}
\end{equation*}
From \eqref{claim} and a change of variables, the left-hand side of \eqref{L2-10} is dominated by
\begin{align*}
  &\sum_{j\geq 0}2^{j/2}\beta_{k,j}\Big\|\eta_{j}(\tau-\omega(\xi))\mathcal{F}_t(\varphi(t)e^{-t\epsilon\xi^2})(\tau-\omega(\xi))\int_\mathbb{R}g_k(\xi,\tau')\,d\tau'
\Big\|_{L^2_{\xi,\tau}}\\
\lesssim &\sum_{j\geq0}2^{j/2}\beta_{k,j}\sup_{\xi\in I_k}\|P_j(\varphi(t)e^{-\epsilon t\xi^2})(t)\|_{L_{t}^2}\cdot \Big\|\int_\mathbb{R}g_k(\xi,\tau')\,d\tau'\Big\|_{L^2_\xi}\\
\lesssim & 2^{k/2}\Big\|\int_\mathbb{R}\chi_{[k-1,k+1]}(\xi)(\tau'-\omega(\xi)+i)^{-1}
\eta_{\leq k}(\tau'-\omega(\xi))\int_\mathbb{R}e^{-ix\xi}h(x,\tau')\,dx\,d\tau'\Big\|_{L^2_\xi}\\
\lesssim & 2^{k/2}\Big\|\chi_{[k-1,k+1]}(\xi)\int_\mathbb{R}e^{-ix\xi}\mathcal{L}_k(h(x,\cdot))(\omega(\xi))dx\Big\|_{L^2_\xi}\\
\lesssim & 2^{k/2}\int_\mathbb{R}\|\chi_{[k-1,k+1]}(\xi)\mathcal{L}_k(h(x,\cdot))(\omega(\xi))\|_{L^2_\xi}dx\lesssim \|h\|_{L^1_xL^2_\tau}\lesssim\|g_k\|_{Y_k},
\end{align*}
the proof of \eqref{L2-10} is completed. Thus we have shown that
\begin{equation*}
\|T\|_{Y_k\to Z_k}\leq C\ \ \ \ \text{ uniformly in }k\geq 100 \ \text{and}\ \epsilon\in [0,1].
\end{equation*}

(2) Case $k=0$.

(2-a) Assume first that $f_0\in X_0$. Similar to $k\geq 1$, we still denote $f_0^\#(\xi,\mu')=f_0(\xi,\mu'+\omega(\xi))$ and $T(f_0)^\#(\xi,\mu)=T(f_0)(\xi,\mu+\omega(\xi))$. Due to $|\xi|\leq 2$, it follows immediately that $\|f_0^\#(\xi,\mu')\|_{X_0}\sim\|f_0\|_{X_0}$ and $\|T(f_0)^\#(\xi,\mu)\|_{X_0}\sim\|T(f_0)\|_{X_0}$. The similar argument as $k\geq 1$, we still divide $T(f_0)^\#(\xi,\mu)$ into four parts:
\begin{align*}
T(f_0)^\#(\xi,\mu)=&\mathcal{F}_t\Big[\varphi(t)\int_{|\mu'|\leq
1}\frac{e^{it\mu'}-1}{i\mu'+\epsilon\xi^2}(\mu'+\omega(\xi)+i)f_0^\#(\xi,\mu')d\mu' \Big](\mu)\nonumber\\
&+\mathcal{F}_t\Big[\varphi(t)\int_{|\mu'|\leq
1}\frac{1-e^{-\epsilon
 t\xi^2}}{i\mu'+\epsilon\xi^2}(\mu'+\omega(\xi)+i)f_0^\#(\xi,\mu')d\mu' \Big](\mu)\nonumber\\
&+\mathcal{F}_t\Big[\varphi(t)\int_{|\mu'|\geq
1}\frac{e^{it\mu'}}{i\mu'+\epsilon\xi^2}(\mu'+\omega(\xi)+i)f_0^\#(\xi,\mu')d\mu' \Big](\mu)\nonumber\\
&-\mathcal{F}_t\Big[\varphi(t)\int_{|\mu'|\geq
1}\frac{e^{-\epsilon
 t\xi^2}}{i\mu'+\epsilon\xi^2}(\mu'+\omega(\xi)+i)f_0^\#(\xi,\mu')d\mu' \Big](\mu)\nonumber\\
:=&I+II+III-IV.
\end{align*}
We first consider the contribution of $I$. By the definition of $X_0$ and Taylor's expansion, we obtain
\begin{align*}
\|I\|_{X_0}\lesssim&
\sum^\infty_{j=0}\sum_{k'\leq 1} \sum_{n\geq
1}\frac{2^{j-k'/2}}{n!}\big\|\eta_j(\mu)\mathcal{F}_t(\varphi(t)t^n)\big\|_{L^2_\mu}\Big\|\int_{|\mu'|\leq
1}\frac{|\mu'|}{|i\mu'+\epsilon\xi^2|}|\chi_{k'}(\xi)f_0^\#(\xi,\mu')|d\mu'\Big\|_{L_{\xi}^2}\nonumber\\
\lesssim&
\sum_{n\geq 1}\frac{\|t^n\varphi(t)\|_{B^{1}_{2,1}}}{n!}\sum^\infty_{j'=0}\sum_{k'\leq 1} 2^{j'/2-k'/2} \|\eta_{j'}(\mu')\chi_{k'}(\xi)f_0^\#(\xi,\mu')\|_{L^2_{\xi,\mu'}}
\lesssim \|f_0^\#(\xi,\mu')\|_{X_0}.
\end{align*}
For $II$, we just take Taylor's expansion to $1-e^{-\epsilon t\xi^2}$, then use the factor $\epsilon\xi^2$ to eliminate the denominator $i\mu'+\epsilon\xi^2$ and get the conclusion similar to $I$.
We then consider the contribution of $III$. Due to the algebraic structure of $B^1_{2,1}$, we know
\begin{align*}
\|III\|_{X_0}\lesssim&\sum^\infty_{j=0}\sum_{k'\leq 1} 2^{j-k'/2}\Big\|\eta_j(\mu)|\mathcal{F}_t\varphi|\star\|\chi_{k'}(\xi)f_0^\#(\xi,\cdot)\|_{L^2_\xi}\Big\|_{L^2_\mu}\nonumber\\
\lesssim&\sum_{k'\leq 1} 2^{-k'/2}\sum^\infty_{j=0}2^j\|P_j(\mathcal{F}^{-1}_\mu|\mathcal{F}_t\varphi|\cdot\mathcal{F}^{-1}_\mu\|\chi_{k'}(\xi)f_0^\#(\xi,\mu)\|_{L^2_\xi})
\|_{L^2_t}\nonumber\\
\lesssim& \sum^\infty_{j=0}\sum_{k'\leq 1} 2^{j-k'/2}\|\eta_j(\mu)\chi_{k'}(\xi)f_0^\#(\xi,\mu)\|_{L^2_{\xi,\mu}}\lesssim \|f_0^\#(\xi,\mu)\|_{X_0}.
\end{align*}
Finally, we consider the contribution of $IV$,
\begin{align*}
\|IV\|_{X_0}\lesssim& \sum^\infty_{j=0}\sum_{k'\leq 1} 2^{j-k'/2}\Big\|\eta_j(\mu)\mathcal{F}_t\Big[\varphi(t)\sum_{n\geq
0}\frac{ t^n(\epsilon\xi^2)^n}{n!}\Big]\int_{|\mu'|\geq
1}|\chi_{k'}(\xi)f_0^\#(\xi,\mu')|d\mu'\Big\|_{L_{\xi,\mu}^2}\nonumber\\
\lesssim&
\sum_{n\geq 0}\frac{C^n\|t^n\varphi(t)\|_{B^{1}_{2,1}}}{n!}\sum_{k'\leq 1} 2^{-k'/2}
\Big\|\int_{|\mu'|\geq1}|\chi_{k'}(\xi)f_0^\#(\xi,\mu')|d\mu'\Big\|_{L_\xi^2}\nonumber\\
\lesssim& \sum^\infty_{j'=0}\sum_{k'\leq 1} 2^{j'/2-k'/2} \|\eta_{j'}(\mu')\chi_{k'}(\xi)f_0^\#(\xi,\mu')\|_{L^2_{\xi,\mu'}}
\lesssim \|f_0^\#(\xi,\mu')\|_{X_0}.
\end{align*}
Now we have obtained that
\begin{equation*}
\|T\|_{X_0\to X_0}\leq C\ \ \ \ \text{ uniformly in } \epsilon\in(0,1].
\end{equation*}

(2-b) Assume now that $f_0=g_{0,j'}(\xi,\tau')\in Y_0$ is supported in $\tilde{I}_0\times I_{j'}$.
We analyze two cases: $j'\geq 5$ and $j'\leq 4$. When $j'\geq 5$, it follows that $|\tau'-\omega(\xi)|\geq 1$, thus the denominator $i(\tau'-\omega(\xi))+\epsilon\xi^2$ is far from origin and there is no singularity.  When $j'\leq 4$, the singularity occurs so that we must handle this case more carefully.

If $j'\geq 5$, we get that $|\tau'-\omega(\xi)|\geq 1$ due to $|\xi|\leq 2$. We rewrite
\begin{align*}
T(g_{0.j'})(\xi,\tau)=&-i\bigg[\int_\mathbb{R}\hat{\varphi}(\tau-\tau')\frac{\tau'+i}{\tau'-\omega(\xi)-i\epsilon\xi^2}
g_{0,j'}(\xi,\tau')d\tau'\nonumber\\
&-\mathcal{F}_t(\varphi(t)e^{-t\epsilon\xi^2})(\tau-\omega(\xi))\int_\mathbb{R}\frac{\tau'+i}{\tau'-\omega(\xi)-i\epsilon\xi^2}g_{0,j'}(\xi,\tau')\,d\tau'\bigg]\nonumber\\
:=&-i(I-II),
\end{align*}
and divide each term into two parts:
\begin{align*}
I=\int_\mathbb{R}\hat{\varphi}(\tau-\tau')\bigg(1+\frac{\omega(\xi)+i\epsilon\xi^2+i}{\tau'-\omega(\xi)-i\epsilon\xi^2}\bigg)
g_{0,j'}(\xi,\tau')d\tau':=I_1+I_2;
\end{align*}
\begin{align*}
II=\mathcal{F}_t(\varphi(t)e^{-t\epsilon\xi^2})(\tau-\omega(\xi))\int_\mathbb{R}\bigg(1+\frac{\omega(\xi)+i\epsilon\xi^2+i}{\tau'-\omega(\xi)-i\epsilon\xi^2}\bigg)
g_{0,j'}(\xi,\tau')d\tau':=II_1+II_2.
\end{align*}

We claim that
\begin{align}\label{0-I1}
\|I_1\|_{Y_0}\lesssim \|g_{0,j'}(\xi,\tau')\|_{Y_0}.
\end{align}
Indeed, by the definition of $Y_0$ and Plancherel's theorem, we have
\begin{align*}
\|I_1\|_{Y_0}&=\sum_{j\geq 1} 2^j \Big\|\eta_j(\tau)\int_\mathbb{R}\hat{\varphi}(\tau-\tau')\mathcal{F}_{\xi}^{-1}[g_{0,j'}(\xi,\tau')]d\tau'\Big\|_{L^1_xL^2_{\tau}}\\
&\quad+\sum_{j\leq 0} \Big\|\chi_{j}(\tau)\int_\mathbb{R}\hat{\varphi}(\tau-\tau')\mathcal{F}_{\xi}^{-1}[(g_{0,j'})(\xi,\tau')]d\tau'\Big\|_{L^1_xL^2_{\tau}}:=I_1^h+I_1^l.
\end{align*}
For $I_1^h$, if $j\geq j'+C$, we know that $|\tau|\sim |\tau-\tau'|$, thus by using Young's inequality, we can get that
 \begin{align*}
I_1^h&\lesssim \sum_{j\geq j'+C} 2^j \Big\|\eta_j(\tau)|\tau|^{-2}\int_\mathbb{R}|\hat{\varphi}(\tau-\tau')(\tau-\tau')^2||\mathcal{F}_{\xi}^{-1}[g_{0,j'}(\xi,\tau')]|d\tau'\Big\|_{L^1_xL^2_{\tau}}\\
&\quad + \sum_{j\leq j'+C} 2^j \Big\|\int_\mathbb{R}\hat{\varphi}(\tau-\tau')\mathcal{F}_{\xi}^{-1}[g_{0,j'}(\xi,\tau')]d\tau'\Big\|_{L^1_xL^2_{\tau}}\\
&\lesssim 2^{j'}\|\mathcal{F}^{-1}(g_{0,j'})\|_{L^1_xL^2_t}=\|g_{0,j'}\|_{Y_0}.
\end{align*}
For $I_1^l$, we could use H\"older's inequality and Young's inequality to obtain that
 \begin{align*}
I_1^l&\lesssim \sum_{j\leq 0} \|\chi_{j}(\tau)\|_{L^2_\tau}\Big\|\int_\mathbb{R}\hat{\varphi}(\tau-\tau')\mathcal{F}_{\xi}^{-1}[(g_{0,j'})(\xi,\tau')]d\tau'\Big\|_{L^1_xL^\infty_{\tau}}\\
&\lesssim \|\mathcal{F}_{\xi}^{-1}[g_{0,j'}(\xi,\tau')]\|_{L^1_xL^2_{\tau'}}\lesssim \|g_{0,j'}\|_{Y_0}.
\end{align*}
Now the claim \eqref{0-I1} is obtained, as desired.

For the term $I_2$, from \eqref{0-I1} we only need to show that
\begin{align}\label{0-I2}
\bigg\|\frac{\omega(\xi)+i\epsilon\xi^2+i}{\tau'-\omega(\xi)-i\epsilon\xi^2}
g_{0,j'}(\xi,\tau')\bigg\|_{Y_0}\lesssim \|g_{0,j'}\|_{Y_0}.
\end{align}
By the definition of $Y_0$, Plancherel's theorem and H\"older's inequalities, to prove \eqref{0-I2}, it suffices to prove that
\begin{align}\label{0-I21}
\bigg\|\int_{\mathbb{R}}e^{ix\xi}\frac{\omega(\xi)+i\epsilon\xi^2+i}{\tau'-\omega(\xi)-i\epsilon\xi^2}\eta_{j'}(\tau')\chi_{\tilde{I}_0}(\xi)d\xi\bigg\|_{L_x^1L_{\tau'}^\infty}\leq C.
\end{align}
Using the facts that $|\tau'-\omega(\xi)|\geq 1$ and $|\xi|\leq 2$, it is easy to get from integration by parts that
\begin{align*}
\bigg|\int_{\mathbb{R}}e^{ix\xi}\frac{\omega(\xi)+i\epsilon\xi^2+i}{\tau'-\omega(\xi)-i\epsilon\xi^2}\eta_{j'}(\tau')\chi_{\tilde{I}_0}(\xi)d\xi\bigg|\lesssim \frac{1}{\langle x\rangle^2},
\end{align*}
which implies \eqref{0-I21}.

The estimates of the term $II$ can be achieved by using the results before. For $II_1$,  from \eqref{k=0-5} we see that
\begin{align}\label{0-I22}
\|II_1\|_{Z_0}\lesssim \Big\|\mathcal{F}^{-1}_\xi\Big[\int_\mathbb{R}
g_{0,j'}(\xi,\tau')d\tau'\Big]\Big\|_{L^1_x}\lesssim2^{j'/2}\|\mathcal{F}^{-1}g_{0,j'}\|_{L^1_xL^2_t}\lesssim\|g_{0,j'}\|_{Y_0}.
\end{align}
Furthermore, \eqref{0-I2} and \eqref{0-I22} lead to $\|II_2\|_{Z_0}\lesssim\|g_{0,j'}\|_{Y_0}$.

If $j'\leq 4$, the singularity occurs by the reason that $i(\tau'-\omega(\xi))+\epsilon\xi^2$ is near origin.   We rewrite
\begin{align*}
T(g_{0,j'})(\xi,\tau)=&\mathcal{F}_t\Big[\varphi(t)\cdot \int_\mathbb{R}\frac{e^{it(\tau'-\omega(\xi))}-1}{i(\tau'-\omega(\xi))+\epsilon\xi^2}
(\tau'+i)g_{0,j'}(\xi,\tau')d\tau' \Big](\tau-\omega(\xi))\nonumber\\
&+\mathcal{F}_t\Big[\varphi(t)\cdot \int_\mathbb{R}\frac{1-e^{-t\epsilon\xi^2}}{i(\tau'-\omega(\xi))+\epsilon\xi^2}
(\tau'+i)g_{0,j'}(\xi,\tau')d\tau' \Big](\tau-\omega(\xi))\nonumber\\
=&\int_\mathbb{R}\frac{\varphi(\tau-\tau')-\varphi(\tau-\omega(\xi))}{i(\tau'-\omega(\xi))+\epsilon\xi^2}
(\tau'+i)g_{0,j'}(\xi,\tau')d\tau'\nonumber\\
&+\mathcal{F}_t\big[\varphi(t)(1-e^{-t\epsilon\xi^2})\big](\tau-\omega(\xi)) \int_\mathbb{R}\frac{\tau'+i}{i(\tau'-\omega(\xi))+\epsilon\xi^2}
g_{0,j'}(\xi,\tau')d\tau' \nonumber\\
:=&A+B.
\end{align*}
Lemma \ref{lemma3} yields that for part $A$ we only need to prove
\begin{align}\label{claim0}
\bigg\|\int_\mathbb{R}\frac{\varphi(\tau-\tau')-\varphi(\tau-\omega(\xi))}{\tau'-\omega(\xi)}
(\tau'+i)g_{0,j'}(\xi,\tau')d\tau'\bigg\|_{Z_0}\lesssim\|g_{0,j'}\|_{Y_0}.
\end{align}
A simple calculation shows that
\begin{align*}
\frac{\varphi(\tau-\tau')-\varphi(\tau-\omega(\xi))}{\tau'-\omega(\xi)}=c\int_0^1\varphi'(\tau-\alpha\tau'-(1-\alpha)\omega(\xi))\,d\alpha.
\end{align*}
Because of $|\xi|\leq 2$ and $|\tau'|\leq C$, we write
\begin{equation*}
\varphi'(\tau-\alpha\tau'-(1-\alpha)\omega(\xi))(\tau'+i)=\varphi'(\tau-\alpha\tau')(\tau'+i)+R(\xi,\tau,\tau'),
\end{equation*}
where
\begin{equation*}
|R(\xi,\tau,\tau')|\leq C\xi^2(1+|\tau|)^{-4}.
\end{equation*}
Therefore, to prove \eqref{claim0}, we just need to show that for any $\alpha \in[0,1]$
\begin{align*}
A_1:=\Big\|\int\varphi'(\tau-\alpha\tau')(\tau'+i)g_{0,j'}(\xi,\tau')d\tau'\Big\|_{Y_0}\lesssim \|g_{0,j'}\|_{Y_0},
\end{align*}
and
\begin{align*}
A_2:=\Big\|\xi^2(1+|\tau|)^{-4}\int |g_{0,j'}(\xi,\tau')|d\tau'\Big\|_{X_0}\lesssim \|g_{0,j'}\|_{Y_0}.
\end{align*}
In fact, for $|\tau|\sim 2^j$ and $|\tau'|\sim 2^{j'}$, if $j\geq 10$, we have $|\tau-\alpha\tau'|\sim |\tau|$. Thus, Minkowski's inequality and H\"older's inequality give that
\begin{align*}
A_1=&\sum_{j\geq 1} 2^j\Big\|\eta_j(\tau)\int\varphi'(\tau-\alpha\tau')(\tau'+i)\mathcal{F}_\xi^{-1}[g_{0,j'}(\xi,\tau')]d\tau'\Big\|_{L^1_xL^2_\tau}\nonumber\\
&+\sum_{j\leq 0} \Big\|\chi_j(\tau)\int\varphi'(\tau-\alpha\tau')(\tau'+i)\mathcal{F}_\xi^{-1}[g_{0,j'}(\xi,\tau')]d\tau'\Big\|_{L^1_xL^2_\tau}\nonumber\\
\lesssim& \Big(\sum_{j\geq 10}2^j\|\eta_j(\tau)(1+|\tau|)^{-2}\|_{L^2_\tau}+1\Big)\|\mathcal{F}_\xi^{-1}[g_{0,j'}(\xi,\tau')]\|_{L^1_{x,\tau'}}\nonumber\\
\lesssim& 2^{j'/2}\|\mathcal{F}^{-1}[g_{0,j'}(\xi,\tau')]\|_{L^1_xL^2_t}\lesssim \|g_{0,j'}\|_{Y_0},
\end{align*}
where we used the fact that $\sum_{j'\leq0}\|\chi_{j'}(\tau)\|_{L^2_\tau}\lesssim 1$. In addition, we can easily get that
\begin{align*}
A_2=&\sum_{j\geq 0}\sum_{k'\leq 1}2^{j-k'/2}\Big\|\eta_j(\tau)\chi_{k'}(\xi)\xi^2(1+|\tau|)^{-4}\int |g_{0,j'}(\xi,\tau')|d\tau'\Big\|_{L^2_{\xi,\tau}}\nonumber\\
\lesssim&\|g_{0,j'}(\xi,\tau')\|_{L^\infty_\xi L^1_{\tau'}}\lesssim 2^{j'/2}\|\mathcal{F}^{-1}[g_{0,j'}(\xi,\tau')]\|_{L^1_xL^2_t}\lesssim \|g_{0,j'}\|_{Y_0}.
\end{align*}
This completes the proof of \eqref{claim0}.

For part $B$, in order to eliminate the singularity,  we will divide it into three sub-parts.  Due to Taylor's expansion, we have
\begin{align*}
1-e^{-t\epsilon\xi^2}=t\epsilon\xi^2-\sum_{n\geq 2}\frac{(-1)^nt^n(\epsilon\xi^2)^{n}}{n!}.\end{align*}
Let
\begin{align*}
B_1:=&\Big(\mathcal{F}_t\big(\varphi(t)(1-e^{-t\epsilon\xi^2})\big)(\tau-\omega(\xi))-\mathcal{F}_t\big(\varphi(t)(1-e^{-t\epsilon\xi^2})\big)(\tau) \Big)\nonumber\\
&\quad\quad\times\int_\mathbb{R}\frac{\tau'+i}{\tau'-\omega(\xi)-i\epsilon\xi^2}
g_{0,j'}(\xi,\tau')d\tau',
\end{align*}
\begin{align*}
B_2:=\sum_{n\geq 2}\frac{\mathcal{F}_t\big(\varphi(t)t^n(\epsilon\xi^2)^n\big)(\tau)}{n!}\int_\mathbb{R}\frac{\tau'+i}{\tau'-\omega(\xi)-i\epsilon\xi^2}
g_{0,j'}(\xi,\tau')d\tau',
\end{align*}
and
\begin{align*}
B_3:=\mathcal{F}_t\big(\varphi(t)t\epsilon\xi^2\big)(\tau)\int_\mathbb{R}\frac{\tau'+i}{\tau'-\omega(\xi)-i\epsilon\xi^2}
g_{0,j'}(\xi,\tau')d\tau'.
\end{align*}
Next we will prove that $B_1,B_2\in X_0$, and $B_3\in Y_0$. For $B_1$, we use $\epsilon\xi^2$, which comes from Taylor's expansion,  to cancel the denominator $\tau'-\omega(\xi)-i\epsilon\xi^2$, and use $\xi^2$, which comes from the mean value theorem,  to absorb the big weight $2^{-k'/2}$ in the definition of $X_0$. Specifically,
by using the mean value theorem and Taylor's expansion, for some $\theta\in[0,1]$, we have
\begin{align*}
\|B_1\|_{X_0}&\lesssim \sum_{j\geq 0}\sum_{k'\leq 1}2^{j-k'/2}\Big\|\eta_j(\tau)\chi_{k'}(\xi)\xi^2\mathcal{F}_t\big(t\varphi(t)
(1-e^{-t\epsilon\xi^2})\big)(\tau-\theta\omega(\xi))\int_\mathbb{R}\frac{(\tau'+i)g_{0,j'}(\xi,\tau')}{\tau'-\omega(\xi)-i\epsilon\xi^2}d\tau'\Big\|_{L^2_{\xi,\tau}}\nonumber\\
&\lesssim \sum_{n\geq 1}\frac{C^n}{n!}\sum_{j\geq 0}2^{j}\|\eta_j(\tau)\mathcal{F}_t\big(\varphi(t)t^{n+1}
\big)(\tau)\|_{L^2_\tau}\Big\|\int_\mathbb{R}\frac{\epsilon\xi^2|g_{0,j'}(\xi,\tau')|}{|\tau'-\omega(\xi)-i\epsilon\xi^2|}d\tau'\Big\|_{L^\infty_{\xi}}\nonumber\\
&\lesssim \sum_{n\geq 1}\frac{C^n\|\varphi(t)t^{n+1}\|_{B^1_{2,1}}}{n!}\cdot2^{j'/2}\|g_{0,j'}(\xi,\tau')\|_{L^\infty_\xi L^2_{\tau'}}\lesssim \|g_{0,j'}\|_{Y_0}.
\end{align*}
For $B_2$, there is a small factor $(\epsilon\xi^2)^2$ as $n\geq 2$. We use one $\epsilon\xi^2$ to cancel the denominator $\tau'-\omega(\xi)-i\epsilon\xi^2$, and another $\epsilon\xi^2$ to absorb the big weight $2^{-k'/2}$ in the definition of $X_0$. Thus we can get
\begin{align*}
\|B_2\|_{X_0}&\lesssim\sum_{n\geq 2} \sum_{j\geq 0}\sum_{k'\leq 1}2^{j-k'/2}\frac{\|\eta_j(\tau)\mathcal{F}_t(\varphi(t)t^n)(\tau)\|_{L^2_\tau}}{n!}\cdot\bigg\|\int_\mathbb{R}\frac{\chi_{k'}(\xi)(\epsilon|\xi|^{2})^n|g_{0,j'}(\xi,\tau')|}
{|\tau'-\omega(\xi)-i\epsilon\xi^2|}d\tau'\bigg\|_{L^2_\xi}\nonumber\\
&\lesssim \sum_{n\geq 2}\frac{C^n\|\varphi(t)t^{n}\|_{B^1_{2,1}}}{n!}\cdot2^{j'/2}\|g_{0,j'}(\xi,\tau')\|_{L^\infty_\xi L^2_{\tau'}}\lesssim \|g_{0,j'}\|_{Y_0}.
\end{align*}
For $B_3$, notice that
$$\frac{\epsilon\xi^2}{\tau'-\omega(\xi)-i\epsilon\xi^2}=i\bigg(1-\frac{\tau'-\omega(\xi)}{\tau'-\omega(\xi)-i\epsilon\xi^2}\bigg),$$
by the proof of Lemma \ref{lemma3}, we have that
$$
\bigg\|\mathcal{F}^{-1}\Big(\frac{\epsilon\xi^2g_{0,j'}(\xi,\tau')}
{\tau'-\omega(\xi)-i\epsilon\xi^2}\Big)\bigg\|_{L^1_xL^2_t}\lesssim \|\mathcal{F}^{-1}\big(g_{0,j'}(\xi,\tau')\big)\|_{L^1_xL^2_t}.$$
Therefore,
\begin{align*}
\|B_3\|_{Y_0}&\lesssim \sum_{j\geq1} 2^j \|\eta_j(\tau)\mathcal{F}_t\big(\varphi(t)t\big)(\tau)\|_{L^2_\tau}\bigg\|\int_\mathbb{R}(\tau'+i)\mathcal{F}^{-1}_\xi\Big(\frac{\epsilon\xi^2g_{0,j'}(\xi,\tau')}
{\tau'-\omega(\xi)-i\epsilon\xi^2}\Big)d\tau'\bigg\|_{L^1_x}\nonumber\\
&\quad+ \sum_{j'\leq 0} \|\chi_{j'}(\tau)\mathcal{F}_t\big(\varphi(t)t\big)(\tau)\|_{L^2_\tau}\bigg\|\int_\mathbb{R}(\tau'+i)\mathcal{F}^{-1}_\xi\Big(\frac{\epsilon\xi^2g_{0,j'}(\xi,\tau')}
{\tau'-\omega(\xi)-i\epsilon\xi^2}\Big)d\tau'\bigg\|_{L^1_x}\nonumber\\
&\lesssim \big(\|\varphi(t)t\|_{B^1_{2,1}}+\|\varphi(t)t\|_{L^1}\big)\cdot 2^{j'/2}\bigg\|\mathcal{F}^{-1}\Big(\frac{\epsilon\xi^2g_{0,j'}(\xi,\tau')}
{\tau'-\omega(\xi)-i\epsilon\xi^2}\Big)\bigg\|_{L^1_xL^2_t}\nonumber\\
&\lesssim 2^{j'/2}\|\mathcal{F}^{-1}\big(g_{0,j'}(\xi,\tau')\big)\|_{L^1_xL^2_t}\lesssim \|g_{0,j'}\|_{Y_0}.
\end{align*}
Now we have proved that
\begin{equation*}
\|T\|_{Y_0\to Z_0}\leq C\ \ \ \ \text{ uniformly in } \epsilon\in(0,1].
\end{equation*}
Therefore, we complete the proof of Lemma \ref{Lemma2}.  $\hfill \Box$

\section{Bilinear estimates}\label{nonlinear}

In this section we state the main bilinear estimates. We show first the dyadic bilinear estimates in spaces $Z_k$.

\begin{lem}\label{Lem1} ($high\times very\ low\rightarrow high$)
Assume $k\geq 20$, $k_2\in[k-2,k+2]$, $f_{k_2}\in Z_{k_2}$, and $f_{0}\in Z_{0}$. Then
\begin{align}\label{bj1}
2^k\big\|\eta_k(\xi)\cdot(\tau-\omega(\xi)+i)^{-1}f_{k_2}\ast f_0\big\|_{Z_k}\lesssim\|f_{k_2}\|_{Z_{k_2}}\|f_{0}\|_{Z_{0}}.
\end{align}
\end{lem}

\begin{lem}\label{Lem2} ($high\times low\rightarrow high$)
Assume $k\geq 20$, $k_2\in[k-2,k+2]$, $f_{k_2}\in Z_{k_2}$, and $f_{k_1}\in Z_{k_1}$ for any $k_1\in[1,k-10]\cap\mathbb{Z}$. Then
\begin{equation}\label{bj2}
2^k\Big\|\eta_k(\xi)(\tau-\omega(\xi)+i)^{-1}f_{k_2}\ast\sum_{k_1=1}^{k-10}f_{k_1}\Big\|_{Z_k}\lesssim
\|f_{k_2}\|_{Z_{k_2}}\sup_{k_1\in[1,k-10]}\|(I-\partial_\tau^2)f_{k_1}\|_{Z_{k_1}}.
\end{equation}
\end{lem}

\begin{lem}\label{Lem3} ($high\times high\rightarrow low$)
Assume $k,k_1,k_2\in\mathbb{Z}_+$, $k_1,k_2\geq k+10$, $|k_1-k_2|\leq 2$, $f_{k_1}\in Z_{k_1}$, and $f_{k_2}\in Z_{k_2}$. Then
\begin{equation}\label{bk2}
\big\|\xi\cdot\eta_k(\xi)\cdot A_k(\xi,\tau)^{-1}f_{k_1}\ast f_{k_2}\big\|_{X_k}\lesssim 2^{-k/4}\|f_{k_1}\|_{Z_{k_1}}\|f_{k_2}\|_{Z_{k_2}}.
\end{equation}
\end{lem}

\begin{lem}\label{Lem4} ($high\times high\rightarrow high\  or\ low\times low\rightarrow low $)
Assume $k,k_1,k_2\in\mathbb{Z}_+$ have the property that
$\mathrm{max}\,(k,k_1,k_2)\leq\mathrm{min}\,(k,k_1,k_2)+30$,
$f_{k_1}\in Z_{k_1}$, and $f_{k_2}\in Z_{k_2}$. Then
\begin{equation}\label{lem2a}
2^k\big\|\eta_k(\xi)\cdot A_k(\xi,\tau)^{-1}f_{k_1}\ast
f_{k_2}\big\|_{Z_k}\lesssim\|f_{k_1}\|_{Z_{k_1}}\|f_{k_2}\|_{Z_{k_2}}.
\end{equation}
Moreover, any spaces $Z_0$ in the right-hand side of \eqref{lem2a}
can be replaced with $\overline{Z}_0$.
\end{lem}

The main proofs of Lemmas \ref{Lem1}-\ref{Lem4} are already given in \cite[Sections 7 and 8]{Io-Ke} and \cite[Lemma 3.3]{Io-Ke2}. The same argument of bilinear estimates as \cite{Io-Ke,Io-Ke2} works, except for the estimates corresponding to $Y_0$. We only need to consider that $Y_0$ appears in the left-hand side of bilinear estimates, since the norm of $Y_0$ in this paper is larger than that in \cite{Io-Ke,Io-Ke2}. Therefore, we only provide a proof of Lemma \ref{Lem4}.

{\bf Proof of Lemma \ref{Lem4}.} We only consider $k=0$ and $k_1,k_2\leq 30$. If $k_1=0$ or $k_2=0$, we may replace the spaces $Z_0$ in the right-hand side of \eqref{lem2a} with $\overline{Z}_0$. A comparison of $X_k$($1\leq k\leq 30$) and $\overline{Z}_0$ indicates that the proofs of the cases $k_1=0$ or $k_2=0$ are identical to the proofs in the corresponding cases $k_1\geq1$ or $k_2\geq 1$. Therefore we may assume $k_1,k_2\geq 1$, $f_{k_1}=f_{k_1,j_1}$ is supported in $D_{k_1,j_1}$, and $f_{k_2}=f_{k_2,j_2}$ is supported in $D_{k_2,j_2}$. Clearly, $\|f_{k_1}\|_{Z_{k_1}}\approx2^{j_1/2}\beta_{k_1,j_1}\|f_{k_1,j_1}\|_{L^2_{\xi,\tau}}\approx
2^{j_1}\|f_{k_1,j_1}\|_{L^2_{\xi,\tau}}$, and
$\|f_{k_2}\|_{Z_{k_2}}\approx2^{j_2/2}\beta_{k_2,j_2}\|f_{k_2,j_2}\|_{L^2_{\xi,\tau}}\approx
2^{j_2}\|f_{k_2,j_2}\|_{L^2_{\xi,\tau}}$.  It suffices to prove that
\begin{align}\label{lem2b}
\big\|\eta_0(\xi)\cdot (\tau+i)^{-1}f_{k_1,j_1}\ast
f_{k_2,j_2}\|_{Y_0}\lesssim 2^{j_1}\|f_{k_1,j_1}\|_{L^2_{\xi,\tau}}\cdot 2^{j_2}\|f_{k_2,j_2}\|_{L^2_{\xi,\tau}}.
\end{align}
By using the definitions, the left-hand side of \eqref{lem2b} is dominated by
\begin{align*}
&\sum_{j\geq 1}2^j\|\mathcal{F}^{-1}[\eta_j(\tau)\eta_0(\xi)\cdot (\tau+i)^{-1}f_{k_1,j_1}\ast
f_{k_2,j_2}]\|_{L^1_xL^2_t}\\
+&\sum_{j'\leq 0}\|\mathcal{F}^{-1}[\chi_{j'}(\tau)\eta_0(\xi)\cdot (\tau+i)^{-1}f_{k_1,j_1}\ast
f_{k_2,j_2}]\|_{L^1_xL^2_t}:=I+II.
\end{align*}
We first estimate the term $I$. If $j\leq 40$, from H\"older's inequality and Plancherel's theorem, we know that
\begin{align*}
I_{j\leq 40}&\lesssim\|\mathcal{F}^{-1}[\eta_j(\tau)\eta_0(\xi)\cdot (\tau+i)^{-1}f_{k_1,j_1}\ast
f_{k_2,j_2}]\|_{L^1_xL^2_t}\\
&\lesssim \|\mathcal{F}^{-1}[f_{k_1,j_1}\ast
f_{k_2,j_2}]\|_{L^1_xL^2_t}\lesssim \|\mathcal{F}^{-1}(f_{k_1,j_1})\|_{L^2_xL^2_t}\|\mathcal{F}^{-1}(f_{k_2,j_2})\|_{L^2_xL^\infty_t}\\
&\lesssim \|f_{k_1,j_1}\|_{L^2_{\xi,\tau}}2^{j_2/2}\|f_{k_2,j_2}\|_{L^2_{\xi,\tau}}.
\end{align*}
If $j\geq 40$, by examining the supports of the functions, we know that $j\leq \max\{j_1,j_2\}+C$. Therefore we assume $j_1=\max\{j_1,j_2\}$ and $j\leq j_1+C$, then
\begin{align*}
I_{j\geq 40}&\lesssim \sum_{j\leq j_1+C}2^j\|\mathcal{F}^{-1}[\eta_j(\tau)\eta_0(\xi)\cdot (\tau+i)^{-1}f_{k_1,j_1}\ast
f_{k_2,j_2}]\|_{L^1_xL^2_t}\\
&\lesssim \sum_{j\leq j_1+C}\|\mathcal{F}^{-1}[f_{k_1,j_1}\ast f_{k_2,j_2}]\|_{L^1_xL^2_t}\lesssim \sum_{j\leq j_1+C} \|f_{k_1,j_1}\|_{L^2_{\xi,\tau}}2^{j_2/2}\|f_{k_2,j_2}\|_{L^2_{\xi,\tau}}\\
&\lesssim 2^{j_1}\|f_{k_1,j_1}\|_{L^2_{\xi,\tau}}\cdot 2^{j_2}\|f_{k_2,j_2}\|_{L^2_{\xi,\tau}}.
\end{align*}
We next estimate the term $II$. By Plancherel's theorem and H\"older's inequality, we achieve that
\begin{align*}
II&\lesssim \sum_{j'\leq 0}\|\chi_{j'}(\tau)(\tau+i)^{-1}\|_{L^2_\tau}\|\mathcal{F}^{-1}[f_{k_1,j_1}\ast
f_{k_2,j_2}]\|_{L^1_xL^1_t}\\
&\lesssim \sum_{j'\leq 0}2^{j'/2}\|\mathcal{F}^{-1}f_{k_1,j_1}\|_{L^2_xL^2_t}\|\mathcal{F}^{-1}f_{k_2,j_2}\|_{L^2_xL^2_t}\\
&\lesssim \|f_{k_1,j_1}\|_{L^2_{\xi,\tau}}\|f_{k_2,j_2}\|_{L^2_{\xi,\tau}}.
\end{align*}
This completes the proof of \eqref{lem2b}. $\hfill\Box$

With these dyadic bilinear estimates in hand, we can use para-product to decompose the bilinear product, and divide it into several cases according to the interactions. The idea is similar to that in \cite[Section 10]{Io-Ke}, so we omit the details and just state the main bilinear estimates for functions in spaces $F^\sigma$.

\begin{prop}\label{bilinear}
If $\sigma\geq 0$ and $u,v\in F^\sigma$ then
\begin{equation}\label{bilinear1}
\|\partial_x(uv)\|_{N^\sigma}\leq C_\sigma(\|u\|_{F^\sigma}\|v\|_{F^0}+
\|u\|_{F^0}\|v\|_{F^\sigma}).
\end{equation}
\end{prop}

\section{Proof of Theorem 1.1}

In this section we complete the proof of Theorem \ref{thm1}. In terms of the uniform estimates Lemma \ref{Lemma1} and Lemma \ref{Lemma2},  bilinear estiamtes Proposition \ref{bilinear}, the proofs of Theorem \ref{thm1} (a), (b) and (c) are similar to that in \cite{Io-Ke2}, thus we only give the ideas.  For any interval $I=[t_0-a,t_0+a]$, $t_0\in\mathbb{R}$, $a\in[0,5/4]$, and $\sigma\geq 0$ we define the normed spaces
\begin{equation*}
\begin{cases}
&F^\sigma(I)=\big\{u\in \mathcal{S}'(\mathbb{R}\times I):||u||_{F^\sigma(I)}:=\inf\limits_{\widetilde{u}\equiv u\text{ on }\mathbb{R}\times I}||\widetilde{u}||_{F^\sigma}<\infty\big\};\\
&N^\sigma(I)=\big\{u\in \mathcal{S}'(\mathbb{R}\times I):||u||_{N^\sigma(I)}:=\inf\limits_{\widetilde{u}\equiv u\text{ on }\mathbb{R}\times I}||\widetilde{u}||_{N^\sigma}<\infty\big\}.
\end{cases}
\end{equation*}
With this notation, the uniform estimates in Lemma \ref{Lemma1} and Lemma \ref{Lemma2} become
\begin{equation}\label{pt1}
\big\|W_\epsilon(t-t_0)\phi\big\|_{F^\sigma(I)}\leq C \|\phi\|_{\widetilde{H}^\sigma},
\end{equation}
and
\begin{equation}\label{pt2}
\Big\|\int_{t_0}^tW_\epsilon(t-s)(u(s))\,ds\Big\|_{F^\sigma(I)}\leq C\|u\|_{N^\sigma(I)}.
\end{equation}
By combining Proposition \ref{bilinear} we obtain
\begin{equation}\label{pt3}
\Big\|\int_{t_0}^tW_\epsilon(t-s)(\partial_x(u\cdot v)(s))\,ds\Big\|_{F^\sigma(I)}\leq C_\sigma(\|u\|_{F^\sigma(I)}\|v\|_{F^0(I)}+\|u\|_{F^0(I)}\|v\|_{F^\sigma(I)})
\end{equation}
for any $u,v\in F^\sigma(I)$, $\sigma\geq 0$. Finally, the estimate \eqref{l5} becomes
\begin{equation}\label{pt4}
\sup_{t\in I}\|u(.,t)\|_{\widetilde{H}^\sigma}\leq C_\sigma \|u\|_{F^\sigma(I)}\ \ \ \ \text{ for any }u\in F^\sigma(I).
\end{equation}

Given $\phi\in B(\delta,\widetilde{H}^0)\cap \widetilde{H}^\infty$, we construct a solution of \eqref{BOB} by iteration:
\begin{equation}\label{pt5}
\begin{cases}
&u^\epsilon_0=W_\epsilon(t)\phi;\\
&u^\epsilon_{n+1}=W_\epsilon(t)\phi-\frac{1}{2}\int_0^tW_\epsilon(t-s)(\partial_x((u^\epsilon_n)^2)(s))\,ds\,\ \ \ \ \text{ for }\,n\in\mathbb{Z}_+.
\end{cases}
\end{equation}
In the following discussion, we assume that $\delta$ is sufficiently small. By using \eqref{pt1} and \eqref{pt3},  we get easily that
\begin{equation}\label{pt6}
\|u^\epsilon_n\|_{F^0([-1,1])}\leq C\|\phi\|_{\widetilde{H}^0}\ \ \ \ \text{ for any }n\in\mathbb{Z}_+,
\end{equation}
by induction over $n$. Using \eqref{pt3} with $\sigma=0$, \eqref{pt5} and \eqref{pt6}, we can show that
\begin{equation}\label{pt7}
\|u^\epsilon_{n}-u^\epsilon_{n-1}\|_{F^0([-1,1])}\leq C2^{-n}\cdot\|\phi\|_{\widetilde{H}^0}\ \ \ \ \text{ for any }n\in\mathbb{Z}_+,
\end{equation}
by induction over $n$. Next, we can obtain that
\begin{align}
&\|u^\epsilon_n\|_{F^{\sigma}([-1,1])}\leq C \|\phi\|_{\widetilde{H}^{\sigma}}\ \ \ \ \ \ \ \ \ \text{ for any }n\in\mathbb{Z}_+\text{ and }\sigma\in[0,2];\label{pt8}\\
&\|u^\epsilon_n\|_{F^{\sigma}([-1,1])}\leq C(\sigma, \|\phi\|_{\widetilde{H}^{\sigma}})\ \ \ \ \text{ for any }n\in\mathbb{Z}_+\text{ and }\sigma\in[0,\infty);\label{pt20}
\end{align}
and
\begin{equation}\label{pt9}
\|u^\epsilon_{n}-u^\epsilon_{n-1}\|_{F^{\sigma}([-1,1])}\leq C(\sigma, \|\phi\|_{\widetilde{H}^\sigma})\cdot 2^{-n}\ \ \ \ \text{ for any }n\in\mathbb{Z}_+\text{ and }\sigma\in[0,\infty).
\end{equation}
For $\sigma\in[0,2]$, the bound \eqref{pt8} and \eqref{pt9} follow in the same way as the bound \eqref{pt6} and \eqref{pt7}, by combining \eqref{pt1}, \eqref{pt3}, and induction over $n$. For $\sigma\geq 2$, we write $\sigma=\sigma_0+\sigma'$, $\sigma'\in\Z_+$, $\sigma_0\in[0,1)$, and argue by induction over $\sigma'$ similar to Section 4 in \cite{Io-Ke2} to complete the proofs of \eqref{pt8} and \eqref{pt9}. Therefore, we can use \eqref{pt9} and \eqref{pt4} to construct
\begin{equation*}
u^\epsilon=\lim_{n\to\infty}u^\epsilon_n\in C([-1,1]:\widetilde{H}^\infty).
\end{equation*}
In view of \eqref{pt5},
\begin{equation*}
u^\epsilon=W_\epsilon(t)\phi-\frac{1}{2} \int_0^tW_\epsilon(t-s)(\partial_x((u^\epsilon)^2(s)))\,ds\ \ \ \ \text{ on }\mathbb{R}\times [-1,1],
\end{equation*}
so $S_\epsilon^\infty(\phi)=u^\epsilon$ is a solution of the initial-value problem \eqref{BOB}, which completes the proof of Theorem \ref{thm1} (a). For Theorem \ref{thm1} (b) and (c), similar to above argument, we can get easily that for $\phi,\phi'\in B(\delta,\widetilde{H}^0)\cap \widetilde{H}^\infty$ then
\begin{equation}\label{pt10}
\sup_{t\in[-1,1]}\|S_\epsilon^\infty(\phi)-S_\epsilon^\infty(\phi')\|_{\widetilde{H}^{\sigma}}\leq C(\sigma, \|\phi\|_{\widetilde{H}^{\sigma}}+\|\phi'\|_{\widetilde{H}^{\sigma}})\cdot \|\phi-\phi'\|_{\widetilde{H}^{\sigma}},
\end{equation}
which implies Theorem \ref{thm1} (b) and (c).

Finally, we prove Theorem \ref{thm1} (d), i.e. the inviscid limit behavior in $\widetilde{H}^\sigma$, $\sigma\geq 0$. Assume $\phi\in B(\delta,\widetilde{H}^0)\cap \widetilde{H}^\sigma$, let $S^\sigma_\epsilon(\phi)$ and $S^\sigma(\phi)$ denote the nonlinear mappings that associate to any initial data $\phi$ the corresponding solutions of the Cauchy problem  \eqref{BOB} and \eqref{BO}. For convenience, we only give the proof of the case $\sigma=0$, since the proofs of the case $\sigma>0$ are similar.   It suffices to prove
\begin{align}\label{pt15}
\lim_{\epsilon\rightarrow0}\|S^0_\epsilon(\phi)-S^0(\phi)\|_{C([-1,1];\widetilde{H}^0)}=0.
\end{align}
We know that
\begin{align}
&u^\epsilon=S^0_\epsilon(\phi)=W(t)\phi- \int_0^tW(t-s)(\partial_x((u^\epsilon)^2(s)/2)-\epsilon\partial_x^2u^\epsilon(s))\,ds;\label{pt11}\\
&u=S^0(\phi)=W(t)\phi-\int_0^tW(t-s)(\partial_x(u^2(s)/2))\,ds, \label{pt12}
\end{align}
where $W(t)\phi=\mathcal{F}_\xi^{-1}e^{it\omega(\xi)}\mathcal{F}_x\phi$ is the solution of the free Benjamin-Ono evolution.
In terms of \eqref{pt2}, \eqref{pt3}, \eqref{pt11}, and \eqref{pt12},  we have
\begin{align}\label{pt13}
&\|u^\epsilon-u\|_{F^0([-1,1])}=\|S^0_\epsilon(\phi)-S^0(\phi)\|_{F^0([-1,1])}\nonumber\\
\lesssim& \big(\|u^\epsilon\|_{F^0([-1,1])}+\|u\|_{F^0([-1,1])}\big)\|u^\epsilon-u\|_{F^0([-1,1])}+\epsilon\|\partial_x^2u^\epsilon\|_{N^0([-1,1])},
\end{align}
Similar to \eqref{pt6}, we have $\|u^\epsilon\|_{F^0([-1,1])}\leq C\delta$, and $\|u\|_{F^0([-1,1])}\leq C\delta$. Combining that with the definitions $N^0$, $F^0$ and \eqref{pt8}, \eqref{pt13} becomes
\begin{align}\label{pt14}
&\|u^\epsilon-u\|_{F^0([-1,1])}\lesssim\epsilon\|\partial_x^2u^\epsilon\|_{N^0([-1,1])}\lesssim \epsilon\|u^\epsilon\|_{F^2([-1,1])}\lesssim \epsilon\|\phi\|_{\widetilde{H}^2}.
\end{align}
In terms of \eqref{pt4}, we have shown that
\begin{align}\label{pt16}
\sup_{t\in[-1,1]}\|S^0_\epsilon(\phi)-S^0(\phi)\|_{\widetilde{H}^0}\leq C\|S^0_\epsilon(\phi)-S^0(\phi)\|_{F^0([-1,1])}\leq C\epsilon\|\phi\|_{\widetilde{H}^2}.
\end{align}
We now prove \eqref{pt15}. $\forall \eta>0$, it  follows from the Lipschitz continuity that there exists a $K>0$ such that
\begin{align}
&\sup_{t\in[-1,1]}\|S^0_\epsilon(P_{\leq K}\phi)-S^0_\epsilon(\phi)\|_{\widetilde{H}^0}\leq \eta/4; \ \ \ \ \forall \epsilon\in(0,1]\\
&\sup_{t\in[-1,1]}\|S^0(P_{\leq K}\phi)-S^0(\phi)\|_{\widetilde{H}^0}\leq \eta/4.
\end{align}
Fixing $K$, by taking $\epsilon=\epsilon(K)$ sufficiently small, we can get from \eqref{pt16} that
\begin{align*}
\sup_{t\in[-1,1]}\|S^0_\epsilon(P_{\leq K} \phi)-S^0(P_{\leq K} \phi)\|_{\widetilde{H}^0}\leq C\epsilon K^2\cdot \|P_{\leq K}\phi\|_{\widetilde{H}^0}\leq \eta/4.
\end{align*}
Therefore, we have
\begin{align*}
\sup_{t\in[-1,1]}\|S^0_\epsilon(\phi)-S^0(\phi)\|_{\widetilde{H}^0}< \eta,
\end{align*}
which implies \eqref{pt15}. The proof of Theorem \ref{thm1} is completed. $\hfill\Box$


\begin{thebibliography}{99}











\bibitem{ABFS}  L. Abdelouhab, J.L. Bona,  M. Felland and J.C. Saut,
Nonlocal models for nonlinear dispersive waves,
 Physica D,  40 (1989) 360-392.

\bibitem{BTB} T.B. Benjamin, Internal waves of permanent form in fluids of great depth, J. Fluid. Mech., 29 (1967) 559-592.

\bibitem{ER} P.M. Edwin and B. Roberts, The Benjamin-Ono-Burgers equation: an application in solar physics, Wave Motion, 8 (1986) 151-158.


\bibitem{Guo-Wang} Z. Guo and B. Wang, Global well-posedness and inviscid limit for the Korteweg-de Vries-Burgers equation, Journal of Differential Equations,  246(10)(2009) 3864-3901.

\bibitem{GLB} Z. Guo, L. Peng, B. Wang and Y. Wang, Uniform well-posedness and inviscid limit for the
Benjamin-Ono-Burgers equation, Advances in Mathematics, 228 (2011) 647-677.


\bibitem{Herr} S. Herr, Well-posedness for equations of Benjamin-Ono type, Illinois Journal of Mathematics,  51(3)(2007) 951-976.


\bibitem{Io-Ke} A.D. Ionescu, C.E. Kenig, Global well-posedness of
the Benjamin-Ono equation in low-regularity spaces, J. Amer. Math.
Soc.,  20(3)(2007) 753-798.

\bibitem{Io-Ke2} A.D. Ionescu, C.E. Kenig, Complex-valued solutions of the Benjamin-Ono equation,
Harmonic analysis, partial differential equations, and related topics, Contemp. Math., Amer. Math. Soc., Providence, RI,  428(2007) 61-74.

\bibitem{Iorio} R.J. Iorio, On the Cauchy problem for the Benjamin-Ono equation, Comm. Partial Differential Equations, 11 (1986) 1031-1081.



\bibitem{KK} C.E. Kenig and K.D. Koenig, On the local well-posedness of the Benjamin-Ono and
modified Benjamin-Ono equations, Math. Res. Lett., 10 (2003)  879-895.

\bibitem{KPV} C.E. Kenig, G. Ponce and L. Vega, On the generalized
Benjamin-Ono equation, Trans. Amer. Math. Soc., 342(1994) 155-172.



\bibitem{Koch0} H. Koch and N. Tzvetkov, On the local well-posedness of the Benjamin-Ono equation in
$H^s(\mathbb{R})$, Int. Math. Res. Not., 2003 (2003) 1449-1464.

\bibitem{Koch} H. Koch and N. Tzvetkov, Nonlinear wave interactions for the Benjamin-Ono equation, Int. Math. Res. Not., 30 (2005) 1833-1847.



\bibitem{molinet2} L. Molinet, A note on the inviscid limit of the Benjamin-Ono-Burgers equation in the energy space, Proceedings of the American Mathematical Society, 141 (8) (2013) 2793-2798.

\bibitem{MRB} L. Molinet and F. Ribaud, On the Cauchy problem for
the generalized Benjamin-Ono equation with small initial data, C. R. Acad. Sci. Paris, Ser. I, 337 (2003)  523-526.

\bibitem{MRB0} L. Molinet and F. Ribaud, Well-posedness results for
the generalized Benjamin-Ono equation with small initial data, J.
Math. Pures Appl., 83 (2004)  277-311.

\bibitem{MRB1} L. Molinet and F. Ribaud, Well-posedness results for
the generalized Benjamin-Ono equation with arbitrary large initial data, Int. Math. Res. Not., 70 (2004)  3757-3795.

\bibitem{MR} L. Molinet, F. Ribaud, On the low regularity of the
Korteweg-de Vries-Burgers equation, Internat. Math. Res. Notices, 37(2002) 1979-2005.

\bibitem{molinet3} L. Molinet, J.C. Saut and N. Tzvetkov, Ill-posedness issues for the Benjamin-Ono and related equations, SIAM J. Math. Anal., 33
(2001) 982-988.

\bibitem{Ono} H. Ono, Algebraic solitary waves in stratified fluids, J. Phys. Soc. Japan, 39 (1975) 1082-1091.

\bibitem{MOT} M. Otani, Bilinear estimates with applications to the generalized Benjamin-Ono-Burgers equations,
Differential Integral Equations,  18(12) (2005) 1397-1426.

\bibitem{Ponce} G. Ponce, On the global well-posedness of the Benjamin-Ono equation, Differential Integral
Equations, 4 (1991) 527-542.




\bibitem{Tao} T. Tao, Global well-posedness of the Benjamin--Ono equation in $H^1(\mathbb{R})$, J. Hyperbolic Differ. Equ.,  1(2004) 27-49.


\bibitem{Vento} S. Vento, Well-posedness and ill-posedness results for dissipative Benjamin--Ono equations, Osaka J. Math.,
48(4)(2011) 933-958.

\end{thebibliography}
\end{document}